\documentclass[3p, number,sort&compress,times,fleqn]{elsarticle}      

\usepackage{amsmath}
\usepackage{amssymb}
\usepackage{mathtools}
\usepackage[colorlinks = true,
linkcolor = blue,
urlcolor  = blue,
citecolor = blue,
anchorcolor = blue]{hyperref}

\usepackage{subfig}
\usepackage{graphicx}
\usepackage[labelfont=bf, justification=justified]{caption}

\usepackage{float}

\usepackage{csml}

\graphicspath{{./figs/}}
\usepackage{url}

\begin{document}
	
\begin{frontmatter}
	
\title{Mollified finite element approximants of arbitrary order and smoothness}
%on polytopic partitions}
%\subtitle{Do you have a subtitle?\\ If so, write it here}

\author[1]{Eky Febrianto}
\author[2]{Michael Ortiz}
\author[1]{Fehmi Cirak\corref{cor1}}
\ead{f.cirak@eng.cam.ac.uk}

\cortext[cor1]{Corresponding author}

\address[1]{Department of Engineering, University of Cambridge, Cambridge, CB2 1PZ, UK }
\address[2]{Graduate Aerospace Laboratories, California Institute of Technology, Pasadena, CA 91125, USA }

\begin{abstract}
The approximation properties of the finite element method can often be substantially improved by choosing smooth high-order basis functions. It is extremely difficult to devise such basis functions for partitions consisting of arbitrarily shaped polytopes. We propose the mollified basis functions of arbitrary order and smoothness for partitions consisting of convex polytopes. On each polytope an independent local polynomial approximant of arbitrary order is assumed. The basis functions are defined as the convolutions of the local approximants with a mollifier. The mollifier is chosen to be smooth, to have a compact support and a unit volume. The approximation properties of the obtained basis functions are governed by the local polynomial approximation order and mollifier smoothness. The convolution integrals are evaluated numerically first by computing the boolean intersection between the mollifier and the polytope and then applying the divergence theorem to reduce the dimension of the integrals. The support of a basis function is given as the Minkowski sum of the respective polytope and the mollifier. The breakpoints of the basis functions, i.e. locations with non-infinite smoothness, are not necessarily aligned with polytope boundaries.  Furthermore, the basis functions are not boundary interpolating so that we apply boundary conditions with the non-symmetric Nitsche method as in immersed/embedded finite elements. The presented numerical examples confirm the optimal convergence of the proposed approximation scheme for Poisson and elasticity problems.  	
\end{abstract}
	
\begin{keyword}
finite elements, polytopic elements, mollifier,  convolution, Voronoi diagrams
\end{keyword}

\end{frontmatter}

%\tableofcontents
\newpage

%
%--------------------------------------------------------------------------------          
\section{Introduction \label{sec:introduction}}
%--------------------------------------------------------------------------------
%
Smooth high-order finite element approximants are often more efficient and, in general, integrate better with prevalent computer-aided geometric design (CAGD) descriptions~\cite{Cirak:2000aa,Hughes:2005aa,Evans:2009aa}. The construction of mesh-based smooth high-order approximants is currently an active area of research as partly motivated by recent academic and industrial interest in isogeometric analysis. In most mesh-based approaches such approximants are defined as the tensor-products of univariate approximants. Such constructions do not generalise to unstructured meshes and auxiliary techniques are needed in the vicinity of the so-called extraordinary vertices where the tensor-product structure breaks. In CAGD a range of ingenious constructions has been conceived to generate smooth high-order approximants around the extraordinary vertices, see the books~\cite{Farin:2002aa,Peters:2008aa} for an overview. Unfortunately, most of these constructions, including~\cite{Scott:2013aa,jia2013reproducing,majeedCirak:2016,toshniwal2017smooth,toshniwal2017multi, zhang2018subdivision,kapl2018construction,zhang2020manifold}, target bivariate manifolds and do not generalise to the arbitrary variate case. Indeed, there are currently no sufficiently flexible and intuitive non-tensor-product arbitrary-variate constructions that can yield smooth polynomial high-order basis functions. By contrast, the proposed mollified approximation scheme over polytopic partitions is easy to construct, is polynomial and can have arbitrary order and smoothness. 

Convolutional techniques are widely used in the analysis and numerics of partial differential equations. The convolution of a function with a mollifier, i.e. a kernel with a unit volume, yields a function that is smoother than the mollifier and the original function. This smoothing property is, for instance, used to recursively define uniform B-splines~\cite{de1986b}, to analyse non-smooth functions and partial differential equations~\cite{hilbert1973mollifier, evans1998partial, adams2003sobolev}, to postprocess finite element solutions~\cite{thomee1977high, bramble1977higher, mirzaee2013smoothness} and to regularise optimisation problems~\cite{sigmund1998numerical, le2011gradient, bletzinger2014consistent}. Indeed, some of the classical meshless methods, like the smoothed particle hydrodynamics (SPH) \cite{lucy1977numerical, gingold1977smoothed} and the  reproducing kernel particle method (RKPM) \cite{liu1995reproducing}, are defined via convolutions, see also the reviews \cite{li2002meshfree, bessa2014meshfree, chen2017meshfree, huerta2018meshfree}. In SPH and RKPM mollifiers are usually referred to as window, weight or influence functions. Different from our mollified approximation scheme, SPH is intrinsically restricted to low order approximants and the RKPM yields high-order and arbitrarily smooth approximants which are rational. The kernels derived in RKPM depend on the local node distribution and are determined so that they can exactly reproduce a polynomial of a given order. Although RKPM was conceived as a meshless method, it is possible to define its mesh-based cousins~\cite{liu2004reproducing} and to blend it consistently with mesh-based B-spline basis functions~\cite{wang2014consistently,valizadeh2015coupled}.
\begin{figure}[b!]
\centering
	\subfloat[][Surface mesh \label{fig:duckMeshCoarse}] {
		\includegraphics[scale=0.09]{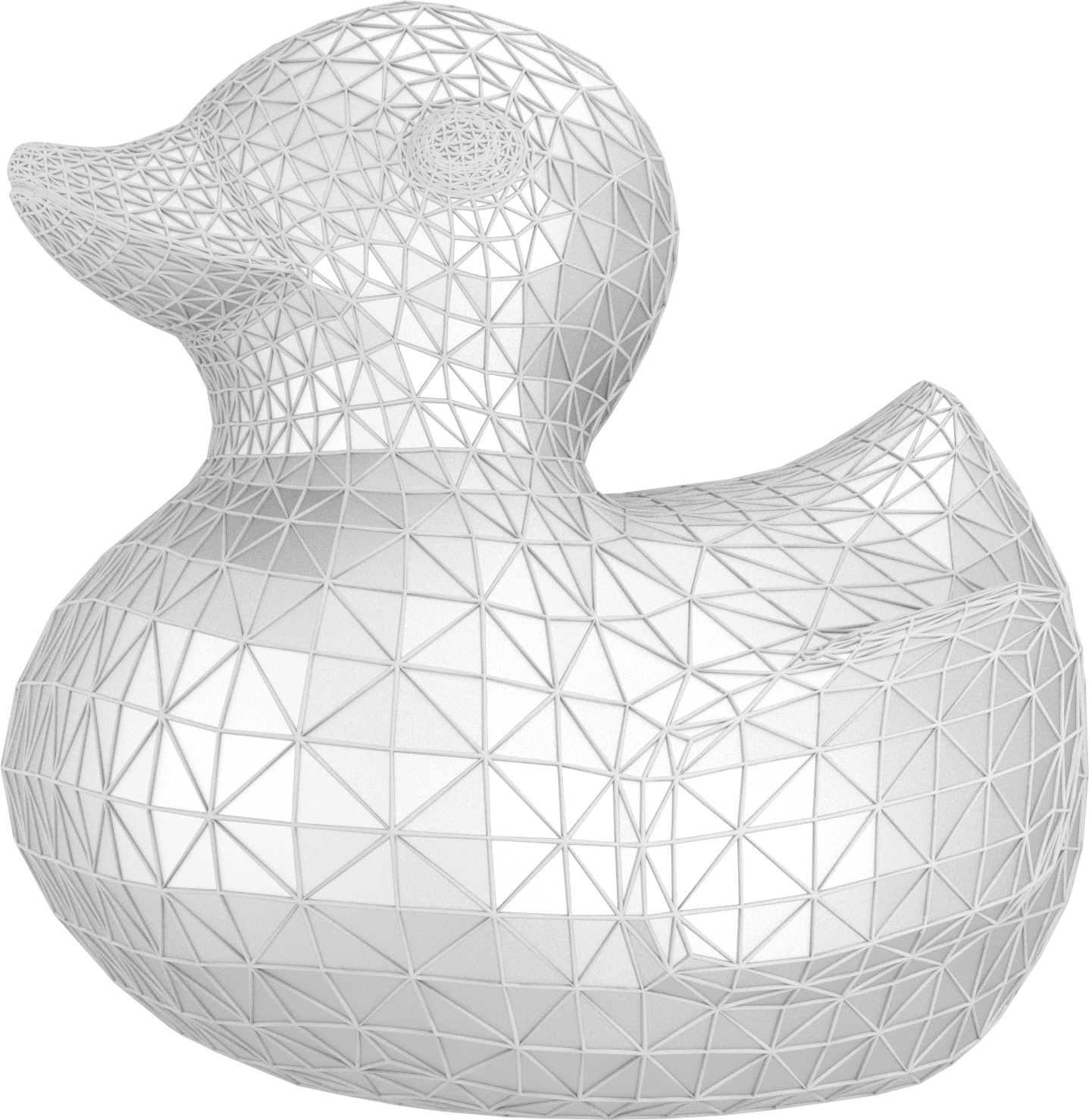} }
\hspace{0.03\textwidth}
	\subfloat[][Clipped Voronoi diagram \label{fig:duckCutoutPerturb}] {
		\includegraphics[scale=0.09]{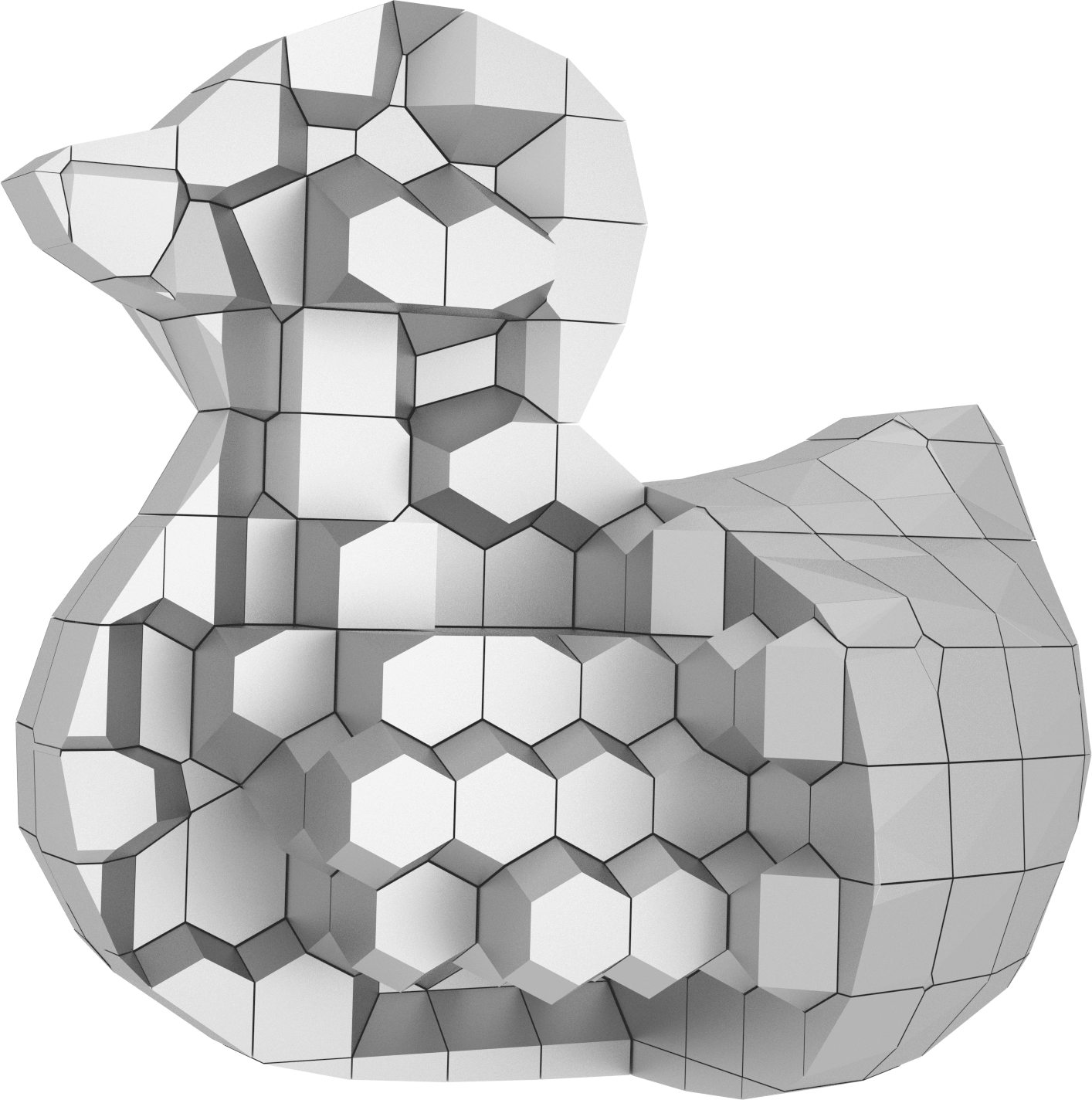}  }
\hspace{0.03\textwidth}
	\subfloat[][Computed finite element potential \label{fig:duckCutoutDispPerturb}] {
		\includegraphics[scale=0.09]{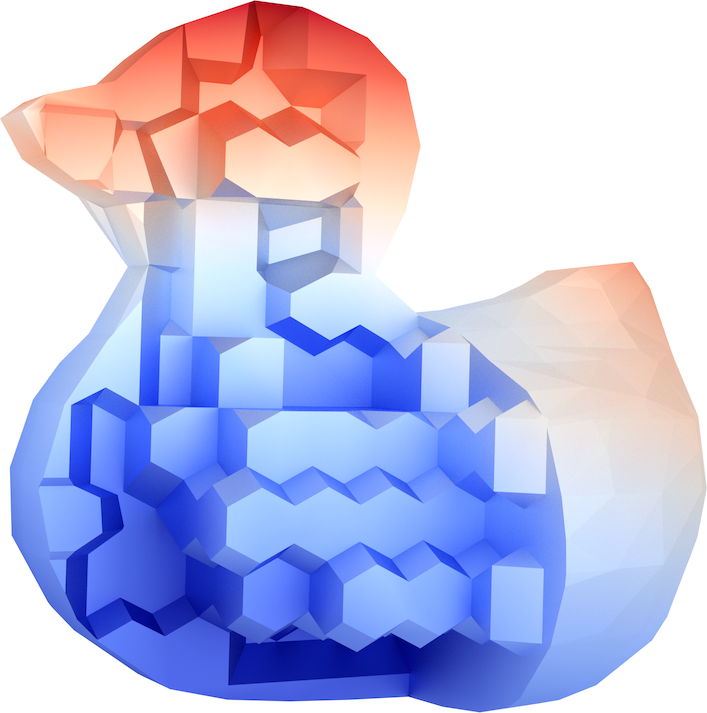}  }
	\caption{Illustrative three-dimensional finite element computation using mollified basis functions. The domain boundary is described with the triangular mesh in (a) and the domain is partitioned with the Voronoi tessellation in (b). The solution of a Poisson problem is shown in (c). Note that in (b) and (c) the cells intersected by the domain boundary have been clipped and others omitted for visualisation purposes.\label{fig:duck}}
\end{figure}

In the proposed mollified approximation scheme each non-overlapping polytopic cell has an independent local polynomial approximant of a prescribed degree~$q^p$. The local approximants are discontinuous across cell boundaries. The convolution of a local approximant with a $C^{k}$ mollifier yields a smoother~$C^{k+1}$ approximant. The chosen mollifiers are compactly supported symmetric polynomials and have a unit volume. It is clear that the convolution of a global polynomial~$f(\vec x)$ of degree~$q^p$ with the chosen mollifiers gives a polynomial~$\widehat{f}(\vec x) \neq f(\vec x)$. However, it is straightforward to find a polynomial~$g(\vec x)$ of degree~$q^p$ such that~$\widehat{g}(\vec x) \equiv f(\vec x)$. This implies that the mollified approximants can exactly reproduce any global function~$f(\vec x)$ of degree~$q^p$. Mollified basis functions for finite element analysis are defined by convolving the local approximant of each cell individually. At a given evaluation point the basis functions are evaluated first by computing the intersection between the support of the mollifier and the cell. Subsequently, the convolution integral over the resulting intersection polytope is evaluated numerically, but exactly (up to round-off error). We apply the divergence theorem to reduce the dimension of the integrals, but other methods of integrating polynomials over polytopes can be used, see e.g.~\cite{chin2015numerical, sudhakar2014accurate}. The obtained basis functions consist of several~$C^{k+1}$ continuously joined polynomial pieces. It is worth emphasising that we use, in contrast to RKPM, a fixed kernel which does not depend on the local node distribution and its convolution with a polynomial~$f(x)$  of degree~$q^p>0$ is not required to yield the same polynomial, i.e.~$\widehat{f}(\vec x) \neq f(\vec x)$. Irrespectively, polynomials $f(x)$ of degree $q^p$ are included in the space spanned by the mollified approximants and can be exactly reproduced.

The derived mollified basis functions can be used as usual in the finite element discretisation of partial differential equations. For ease and efficiency of implementation we assume that each of the polytopic cells representing a finite element is convex. The required geometric operations, like the intersection computations, are significantly simplified by the convexity assumption. We partition the problem domain into a set of convex polytopic cells using a Voronoi diagram, see Figure~\ref{fig:duck}. The evaluation of the finite element integrals requires some care because the~$C^{k+1}$ continuous breakpoints of the basis functions are not aligned with cell boundaries. We use to this end the variationally consistent integration approach proposed in~\cite{chen2017meshfree}, which significantly reduces the number of needed integration points. The required support of the mollified basis functions is the Minkowski sum of the mollifier support with the respective cell~\cite{deBerg2010}. Furthermore, the present version of the mollified basis functions is non-boundary-interpolating so that the Dirichlet boundary conditions are applied weakly with the non-symmetric Nitsche method~\cite{ Nitsche:1971aa, oden1998discontinuoushpfinite, boiveau2015penalty, schillinger2016non} as in immersed/embedded finite element methods, see e.g.~\cite{ Ruberg:2011aa}.

The outline of this paper is as follows. In Section~\ref{sec:preliminaries} we briefly review the convolution of univariate polynomials with a mollifier and characterise the properties of the resulting mollified polynomials.  Subsequently we derive in Sections~\ref{sec:basisUni} and~\ref{sec:basisMulti} first the univariate and then the multivariate mollified basis functions. The key difference between the two cases lies in the evaluation of the convolution integrals. In the univariate case the integrals are evaluated analytically and in the multivariate case numerically. The use of the derived mollified basis functions in finite element analysis, especially the integration and treatment of boundary conditions, is discussed in Section~\ref{sec:fem}. Finally, in Section~\ref{sec:examples} we introduce several Poisson and elasticity examples to confirm the optimal convergence of the developed approach. The paper is supplemented by four appendices which provide convergence estimates and discuss implementation details. 

%
%--------------------------------------------------------------------------------          
\section{Preliminaries \label{sec:preliminaries}}
\begin{figure}[tb]
\centering
	\subfloat[][Constant $f(x)$ with $q^p=0$ \label{fig:mollifiedConst}] {
		\includegraphics[scale=0.38]{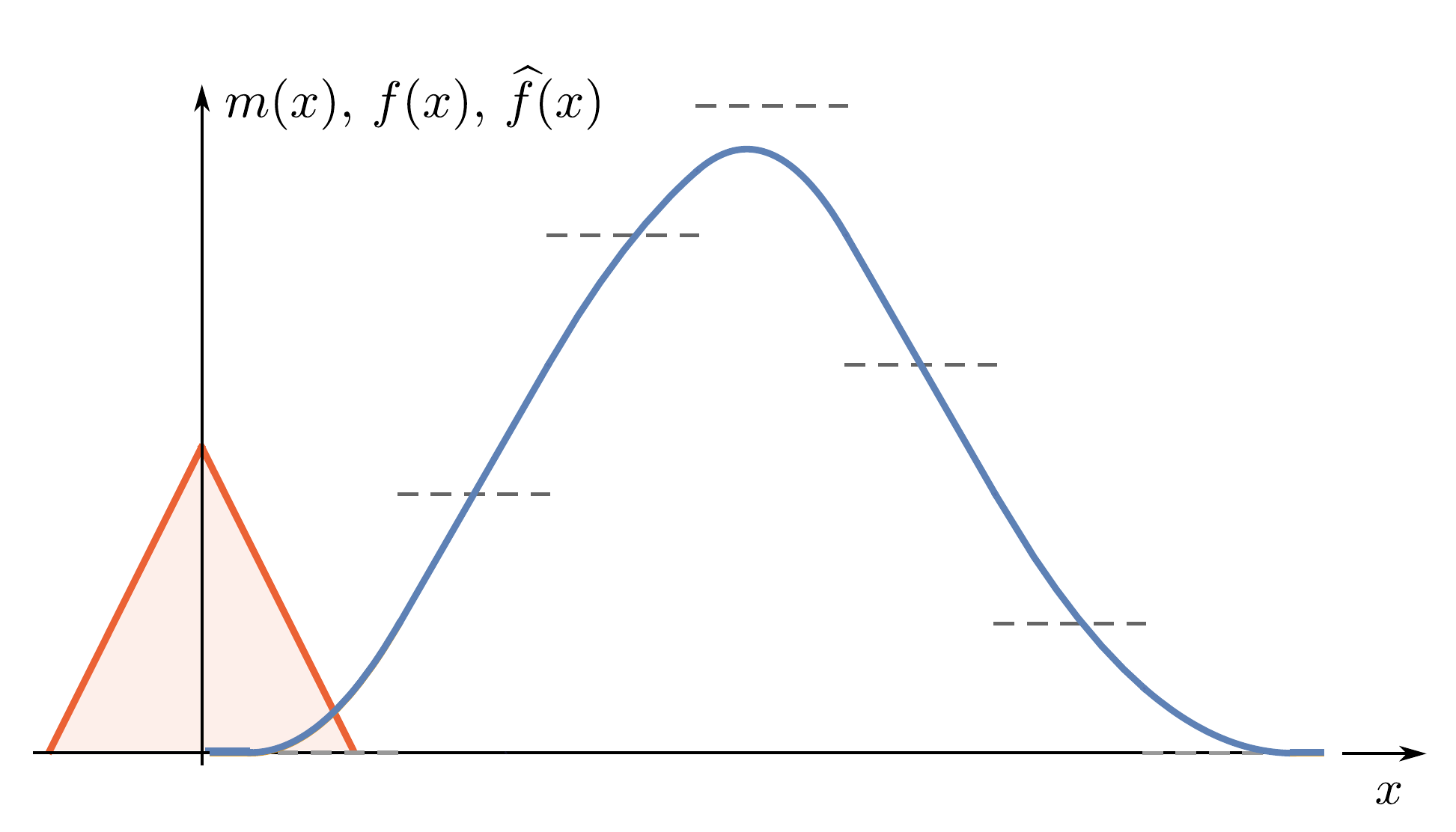} }
\hspace{0.03\textwidth}
	\subfloat[][Linear~$f(x)$ with $q^p=1$ \label{fig:mollifiedLin}] {
		\includegraphics[scale=0.38]{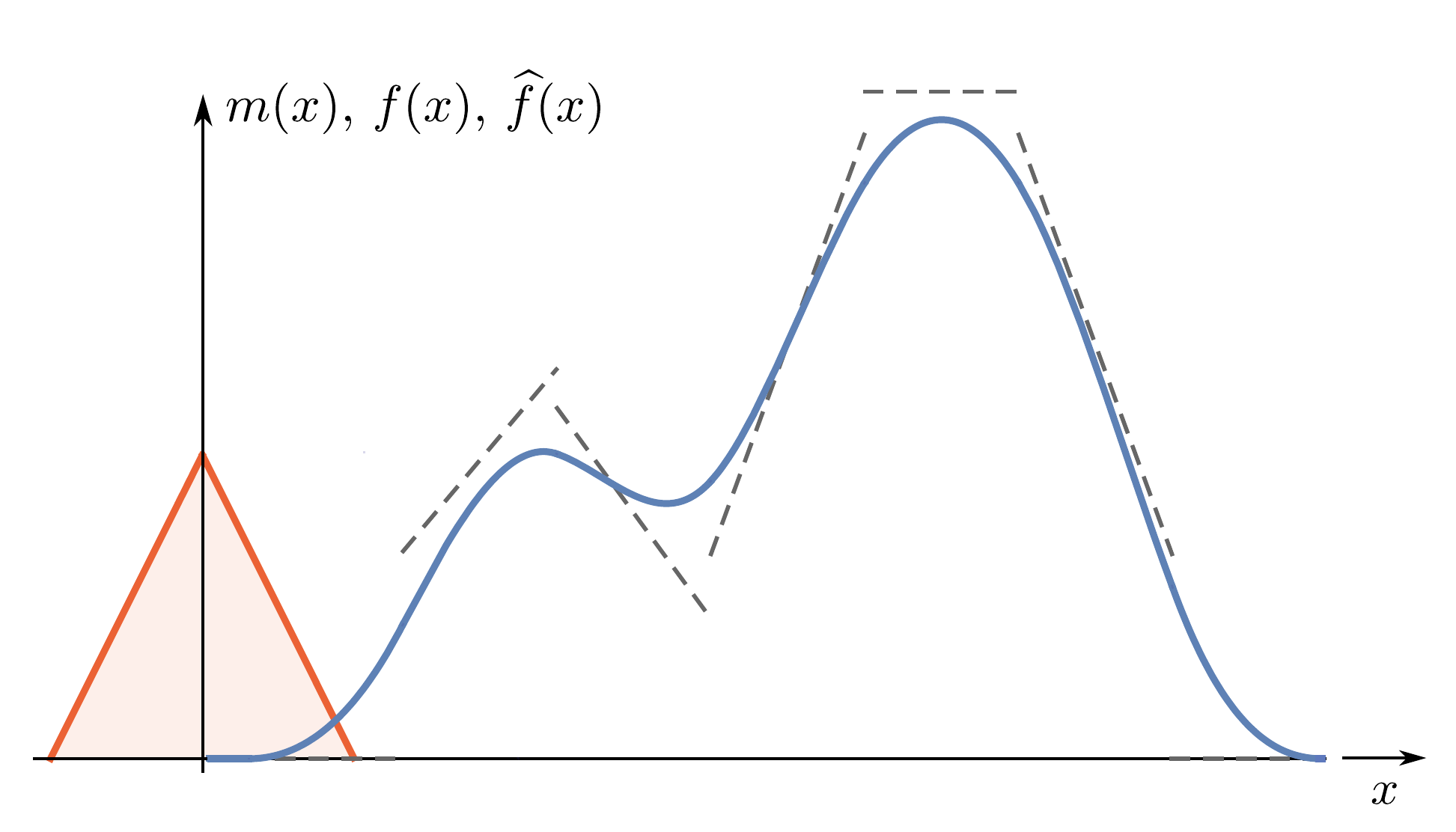}  }
	\caption{Mollification of piecewise discontinuous functions $f(x)$ (black, dashed) with a linear mollifier $m(x)$ (red, solid).  The resulting mollified functions~$\widehat f(x)$ (blue, solid)  are $C^1$ continuous. \label{fig:mollified}}
\end{figure}

%--------------------------------------------------------------------------------
%
We consider the one-dimensional domain~$\Omega \in \mathbb{R}^1$ partitioned into a set of~$n_c$ non-overlapping segments~$\{ \Omega_i \}$, referred to as cells, such that
\begin{equation}
	\Omega = \bigcup_{i=1}^{n_{c}} \Omega_i  \, .
\end{equation}
On each cell~$\Omega_i$  a compactly supported local polynomial is defined,  
\begin{equation} \label{eq:fi}
	f_i(x) = 
	\begin{cases}
		\vec p_i (x)  \cdot \vec \alpha_i \quad  & \text{if } \,  x \in \Omega_i \\
		0						 &  \text{if } \,  x \not \in \Omega_i 
	\end{cases} \, , 
\end{equation}
where the vector~$\vec p_i(x)$ represents a polynomial basis of degree~$q^p$ and~$\vec \alpha_i$ are its polynomial coefficients. See Figure~\ref{fig:mollified} for an illustrative example. We choose in each cell the same polynomial basis~$\vec p_i(x)$, although it is possible to change the type of basis and its polynomial degree.  The sum of the local polynomials defined over the entire domain~$\Omega$ is given by 
\begin{equation} \label{eq:f}
	f(x) = \sum_i \vec p_i(x) \cdot \vec \alpha_i \, . 
\end{equation}
Evidently, across the cell boundaries this function can be discontinuous, i.e.~$f(x) \in C^{-1}$. 
 
The smoothness of~$f(x)$ is increased by convolving it with a mollifier~$m(x)$. The mollifier is chosen such that it has the following properties
\begin{subequations}
\begin{align}
	m(x) & \ge 0  \quad \forall x \in \Omega \\ 
	  \supp m(x)  &=  ( -h_m/2, \, h_m/2 )    \\ 
	\int_\Omega m(x) \, \D x& = 1 \, .
\end{align}
\end{subequations}
That is, the mollifier is non-negative, has a unit volume and a finite support of size~$h_m$. In addition we require that the mollifier is symmetric, i.e.~$m(x) = m(-x)$, and that it has certain smoothness properties, as yet to be specified.  The mollification of  $f(x)$ is defined with the convolution
\begin{equation} \label{eq:mollify}
	\widehat f (x) = m(x) * f(x)  = \int_\Omega m(x-y) f (y) \D y  =   \int_\Omega m(y) f (x-y) \D y  \, .
\end{equation}
The equality of both integrals can be shown by a simple substitution. We usually use the first integral expression in the following. Furthermore, we choose polynomial mollifiers~$m(x)$ of degree~$q^m$ and~$f(x)$ is, as stated above, of degree~$q^p$. The mollified function~$\widehat{f}(x)$ has monomials up to degree~$q^m+q^p+1$. 

If the derivative of the mollifier~$m(x)$ exist, the derivative of the mollified function~$\widehat{f}(x)$ is given by
\begin{equation} \label{eq:mollifyDeriv}
	\frac{\D }{ \D x} \widehat{f}(x) =  \int_\Omega \frac{\D m(x-y)}{\D x} f (y) \D y \, .
\end{equation}
The higher order derivatives are computed similarly. Considering that~${f}(x) \in C^{-1}$ is discontinuous and the mollifier is~$m(x) \in C^{k}$ we can deduce for the mollified smooth function~$\widehat{f}(x) \in C^{k+1}$ .

Finally, the cell-wise definition of~$f(x)$ introduced in~\eqref{eq:mollify} and~\eqref{eq:mollifyDeriv} yields
\begin{equation} \label{eq:mollifyCells}
	\widehat{f}(x) = \sum_i \vec \alpha_i \cdot \int_{\Omega_i} m(x-y) \vec p_i(y) \D y 
\end{equation}
and 
\begin{equation}
	\frac{\D }{ \D x} \widehat{f}(x)  = \sum_i \vec \alpha_i \cdot \int_{\Omega_i}   \frac{\D m(x-y)}{\D x}  \vec p_i(y) \D y  \, .
\end{equation}
The function~$\widehat{f} (x)$ is composed of infinitely smooth polynomial pieces that are smoothly connected, i.e.~$C^{k+1}$, across a finite number of breakpoints. In general the location of the breakpoints does not coincide with the cell boundaries. However, as is known from B-splines, on uniformly partitioned domains, it is possible to choose the support size of the mollifier such that the breakpoints fall on the cell boundaries~\cite{de1986b, Sabin:2010aa}. In the illustrative example in Figure~\ref{fig:mollified} the influence of the choice of the local polynomials~$f_i(x)$ on the mollified function~$\widehat{f} (x)$ is demonstrated. In this example the convolution integral~\eqref{eq:mollifyCells} has been evaluated analytically.

A final remark concerns the reproduction of polynomials with mollified functions. It can be shown that it is possible to find for a given polynomial~$f(x)$ of degree~$q^p$ a polynomial~$g(x)$ of the same degree which yields after mollification~\mbox{$\widehat{g}(x) = m(x) * g(x) \equiv  f(x)$.} Specifically, the mollification of a polynomial 
\begin{equation}
	g(x) = \alpha_0 + \alpha_1 x + \alpha_2 x^2 + \alpha_3 x^3 + \dotsc 
\end{equation}
is given by
\begin{equation} \label{eq:mollifiedPoly}
	\widehat{g}(x) = \alpha_0 + \alpha_1 x + \alpha_2 (x^2 + m_2)  + \alpha_3 ( x^3 + 3 x m_2) + \dots \, ,
\end{equation}
where~$m_s$ are the moments of the mollifier defined as
\begin{equation}
	m_s = \int_\Omega m(x) x^s \D x \, .
\end{equation}
After rearranging the terms in~\eqref{eq:mollifiedPoly} according to the powers of~$x$ it is easy to see how to choose a function~$g(x)$ which is after mollification equal to the given function~$ f(x)$. This implies that the polynomials of degree~$q^p$ are included in the space spanned after mollification.

%
%--------------------------------------------------------------------------------          
\section{Smooth piecewise basis functions \label{sec:basis}}
%--------------------------------------------------------------------------------
%
We now use the mollification approach to derive basis functions on one- and multi-dimensional domains. Again the domain~$\Omega \in \mathbb R^d$ is assumed to be partitioned into a set of  non-overlapping convex polytopes~$\{\Omega_i \}$,  which are in the present paper obtained from a Voronoi diagram. In the following we refer to the polytopes as cells. The mollification approach yields a set of basis functions for each cell. The convolution integrals for obtaining the basis functions are evaluated analytically in the one-dimensional, i.e. univariate, case with $d=1$ and numerically in the multi-dimensional, i.e. multivariate, case with $d \ge 2$. Convergence estimates for the obtained mollified basis are provided in~\ref{app:convergence}.
%
%--------------------------------------------------------------------------------          
\subsection{Univariate basis functions  \label{sec:basisUni}}
%--------------------------------------------------------------------------------
%
To set the stage for multivariate basis functions, we first consider the derivation of univariate basis functions.  The mollified basis functions belonging to a cell~$\Omega_i$ are defined according to the mollification~\eqref{eq:mollifyCells} by 
\begin{equation} 
	f^m(x) = \sum_i \vec \alpha_i \cdot \vec N_i(x) \, ,
\end{equation}
where the vector of mollified basis functions $\vec N_i(x)$ is given by
\begin{equation} \label{eq:basis1d}
	 \vec N_i(x) =  \int_{\Omega_i} m(x-y) \vec p_i (y)  \D y  \, .
\end{equation}
Note that the local polynomial basis is outside the cell~$\Omega_i$ zero and the integration domain is restricted to~$\Omega_i$. As a local polynomial basis~$\vec p_i(y)$ different basis choices are possible, such as the monomial, Lagrange or Bernstein. While this choice has no influence on the approximation quality of the resulting mollified basis, it affects the interpretability of the coefficients~$\vec \alpha_i$ and the conditioning of the resulting finite element system matrices. In our examples we use in each cell~$\Omega_i$ the scaled and shifted monomial basis 
\begin{equation}
	\vec p_i(x) = \begin{pmatrix}
		1 & x & x^2  & \dotsc & x^{q^p} 
		\end{pmatrix}
	\quad \text{  with } x =  \frac{2(x-c_i)}{h} \, ,
\end{equation}
where~$q^p$ is its degree,~$c_i$ is the centre of the cell and~$h$ is the average length of all the cells~$\{ \Omega_i\}$ in the domain. The scaling by $2/h$ ensures that the obtained mollified basis functions have a similar maximum value. It is possible to apply a different scaling factor and to choose the scaling in each cell differently. 

\begin{figure}[]
\centering
	\subfloat[][Mollifier (narrow) and local monomial basis] {
		\includegraphics[width=0.4 \textwidth]{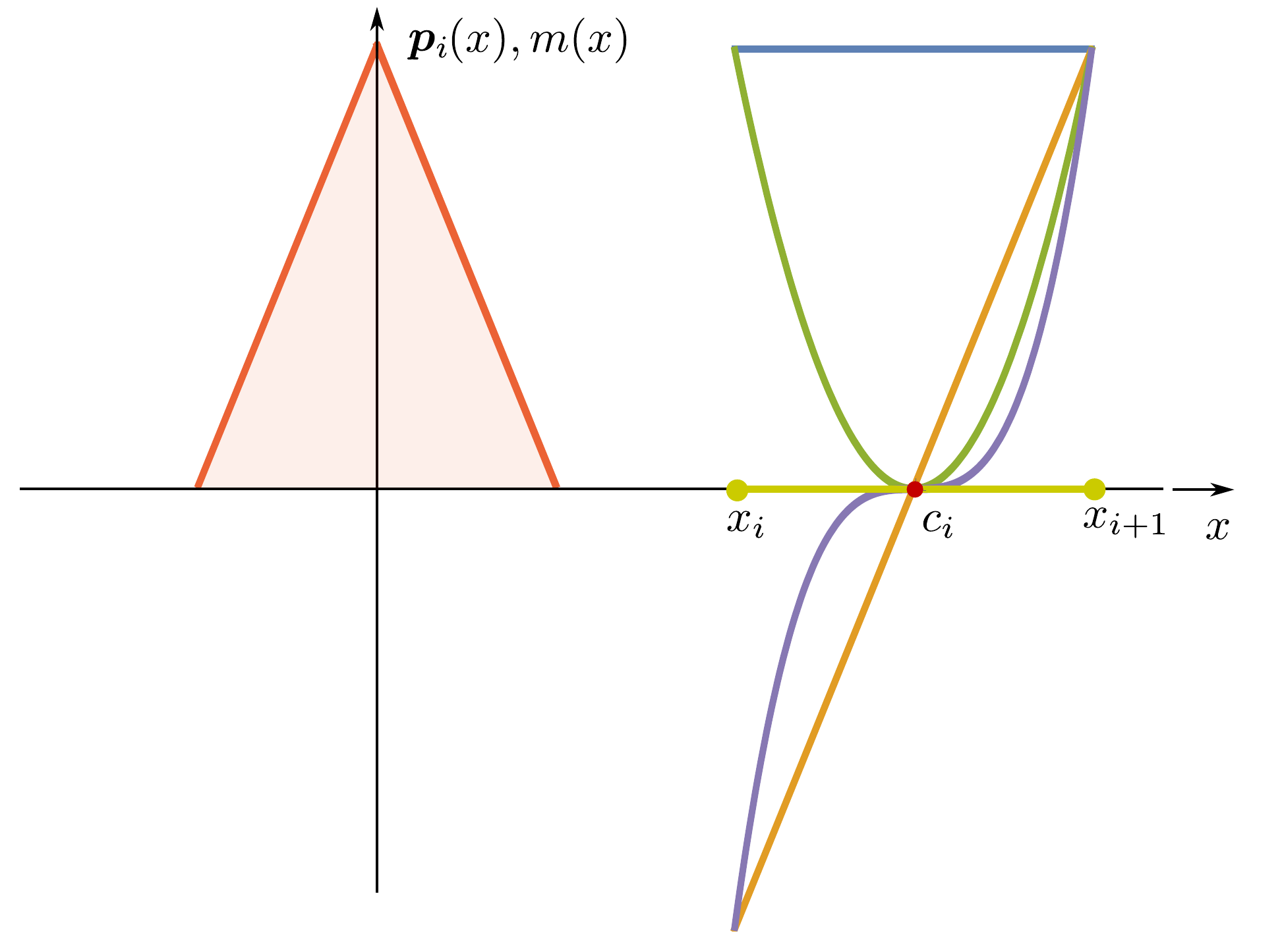} } 
		\hfill
	\subfloat[][Basis functions obtained with the mollifier in (a)] {
		\includegraphics[width=0.4 \textwidth]{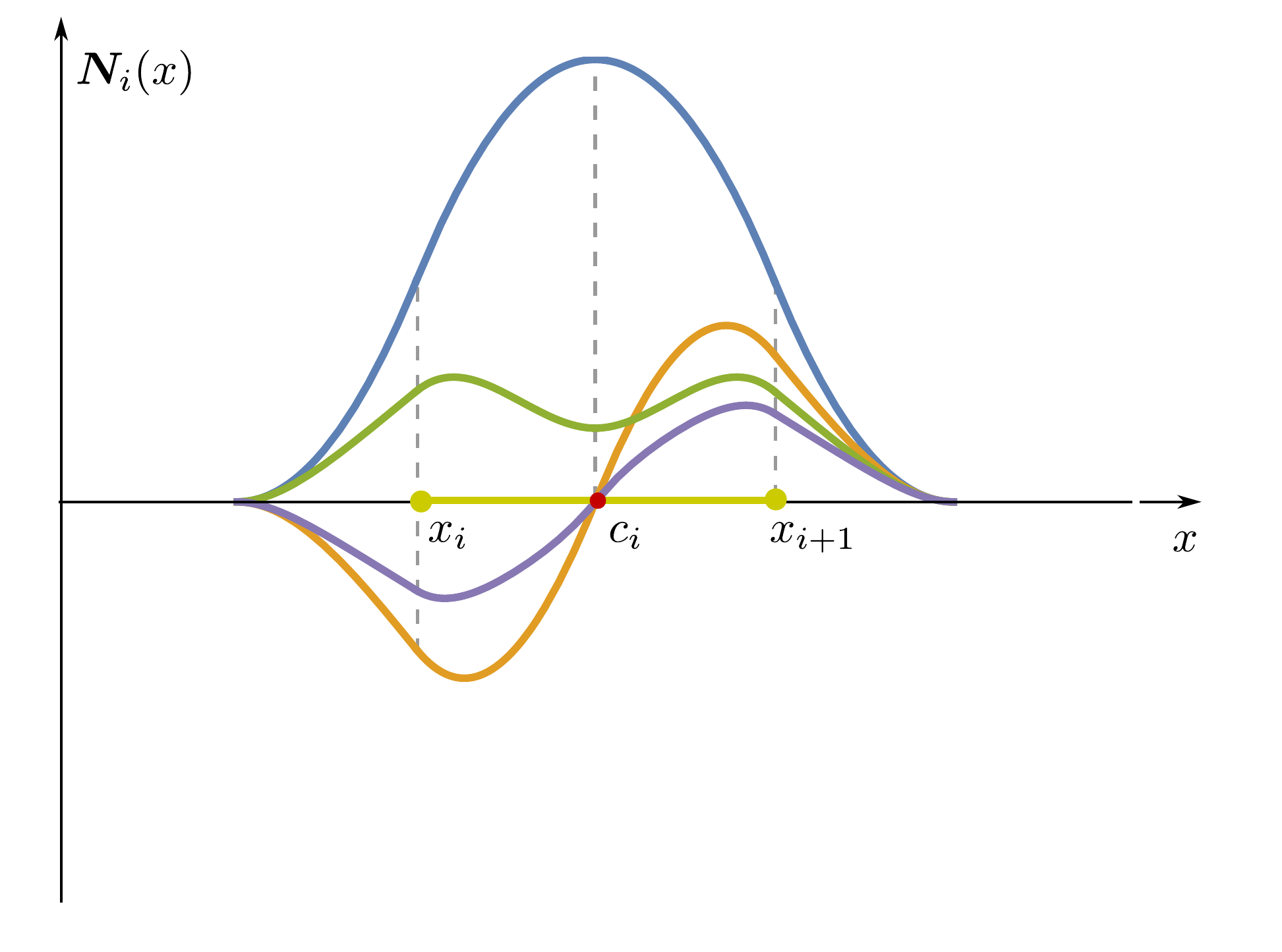} } \\
	\subfloat[][Mollifier (wide) and local monomial basis] {
		\includegraphics[width=0.4 \textwidth]{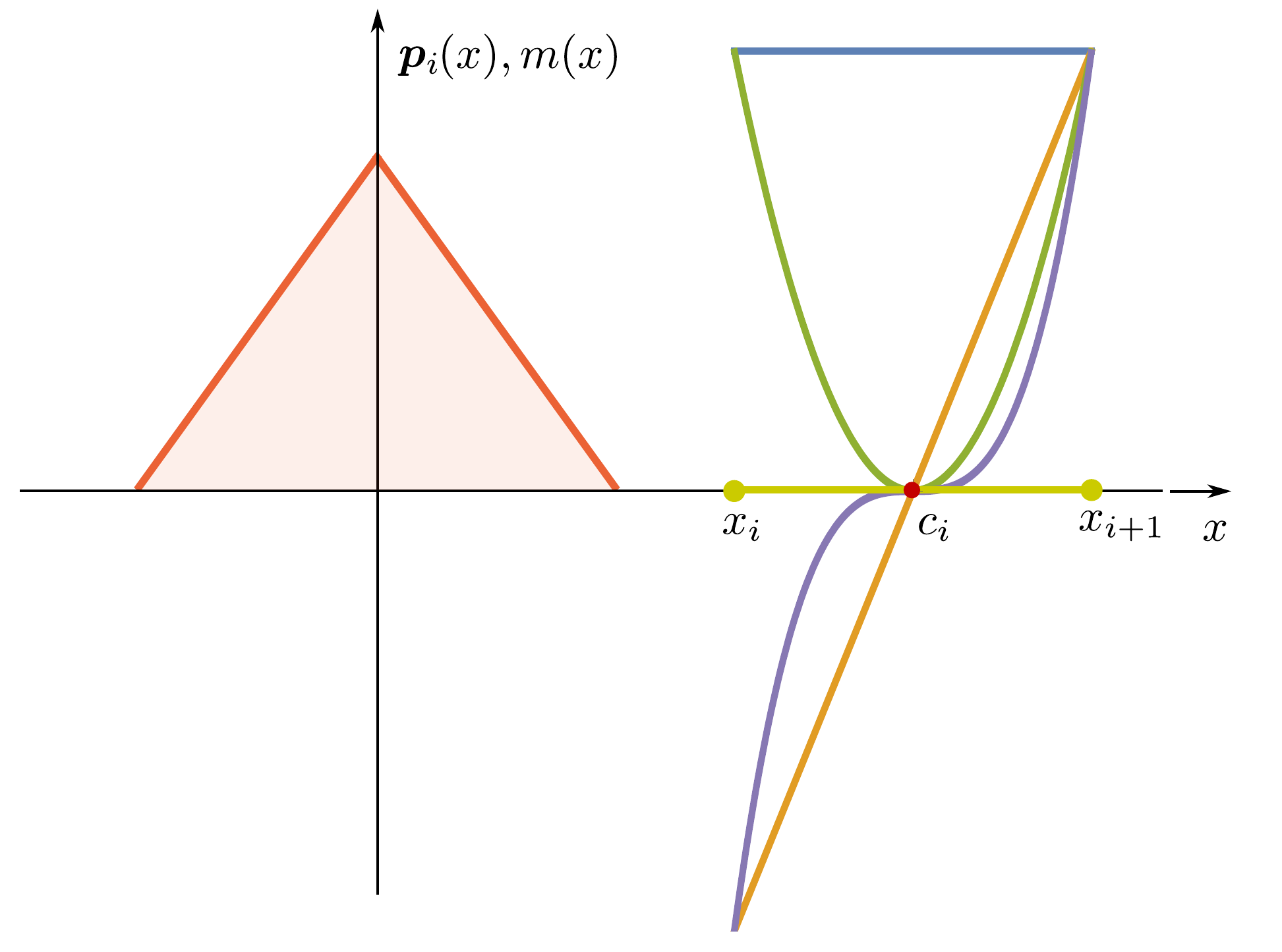} } 
		\hfill
	\subfloat[][Basis functions obtained with the mollifier in (c)] {
		\includegraphics[width=0.4 \textwidth]{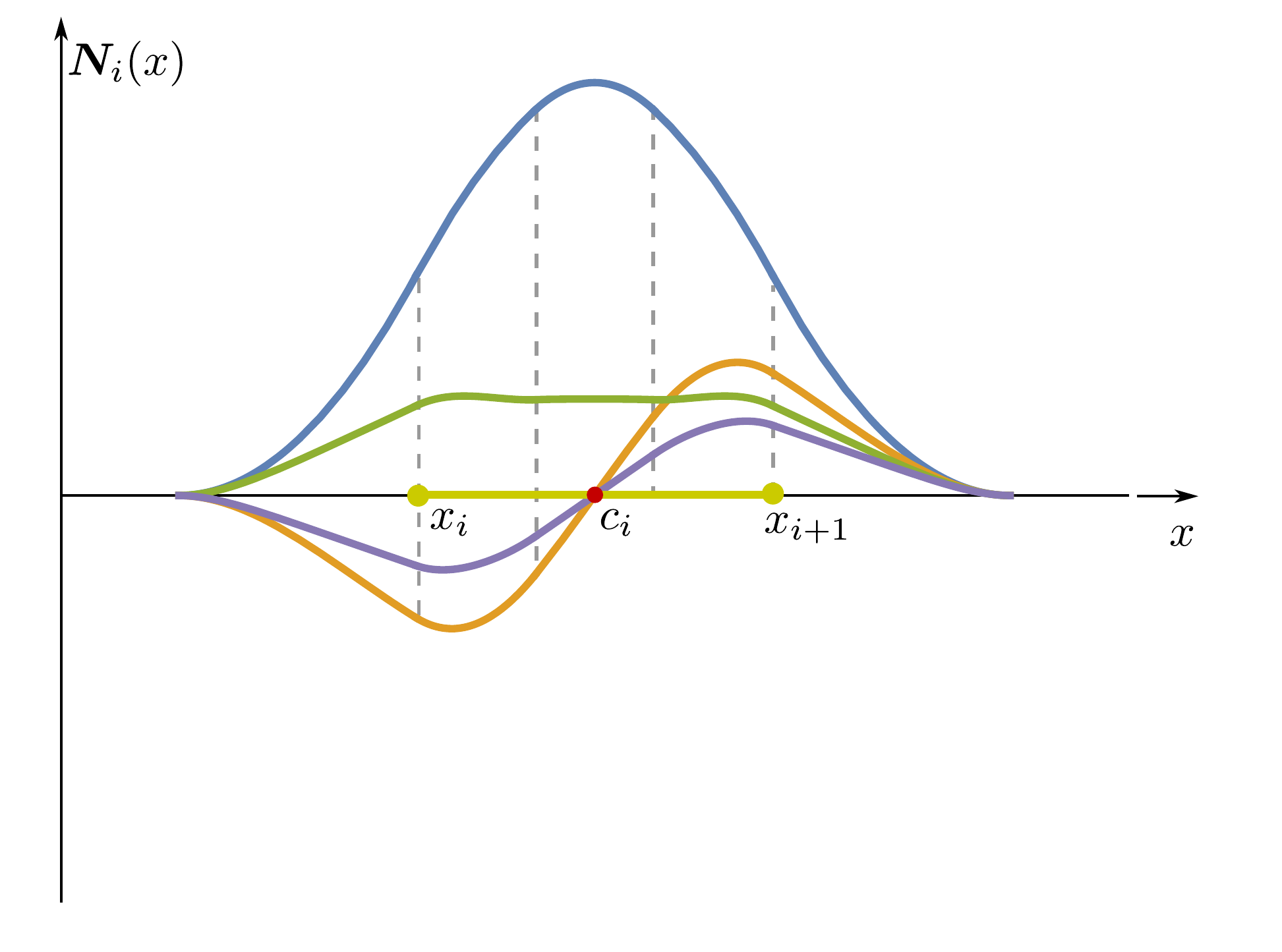} }
	\caption{Univariate mollified basis functions with a bilinear mollifier ($q^m =1$) and a cubic scaled and shifted monomial basis ($q^p=3$) on a cell~$\Omega_i=(x_{i+1}, \, x_{i})$. The mollified basis functions in~(b) are obtained with the narrow mollifier in~(a) and the ones in~(d) with  the wide mollifier in~(c).  The dashed lines in~(b) and~(d) indicate the $C^1$ continuous breakpoints of the obtained basis functions. \label{fig:basis1D}}
\end{figure}

For any given point~$x \in \Omega$ in the domain the mollified basis functions~$\vec N_i(x)$ are evaluated by computing the convolution integral~\eqref{eq:basis1d}. Evidently, when the chosen local polynomial basis~$\vec p_i(x)$ is a monomial basis the mollified basis functions are simply the moments of the mollifier.  The derivatives of the mollified basis functions are computed according to~\eqref{eq:mollifyDeriv}.  In Figure~\ref{fig:basis1D} the mollified basis functions for a cubic monomial basis with~$q^p=3$ and a piecewise linear mollifier~$q^m=1$ with two different support sizes~$h_m$ are shown. In each case the support size of the mollified basis functions is~$h_m+h_{c, i}$, with  the mollifier size~$h_m$ and the cell size $h_{c, i}= x_{i+1} -x_i$. The basis functions consist of several polynomial pieces that are $C^1$ continuously connected at the breakpoints.  The number of breakpoints in each cell depends on the number and arrangement of the breakpoints in the mollifier. Note, although not shown in the figure, the basis functions of the neighbouring cells are non-zero in the considered cell as well. The breakpoints of those neighbouring basis functions may not coincide with the breakpoints of the shown basis functions.

%
%--------------------------------------------------------------------------------          
\subsection{Multivariate basis functions  \label{sec:basisMulti}}
\begin{figure} 
	\centering
	\includegraphics[width=0.4\textwidth]{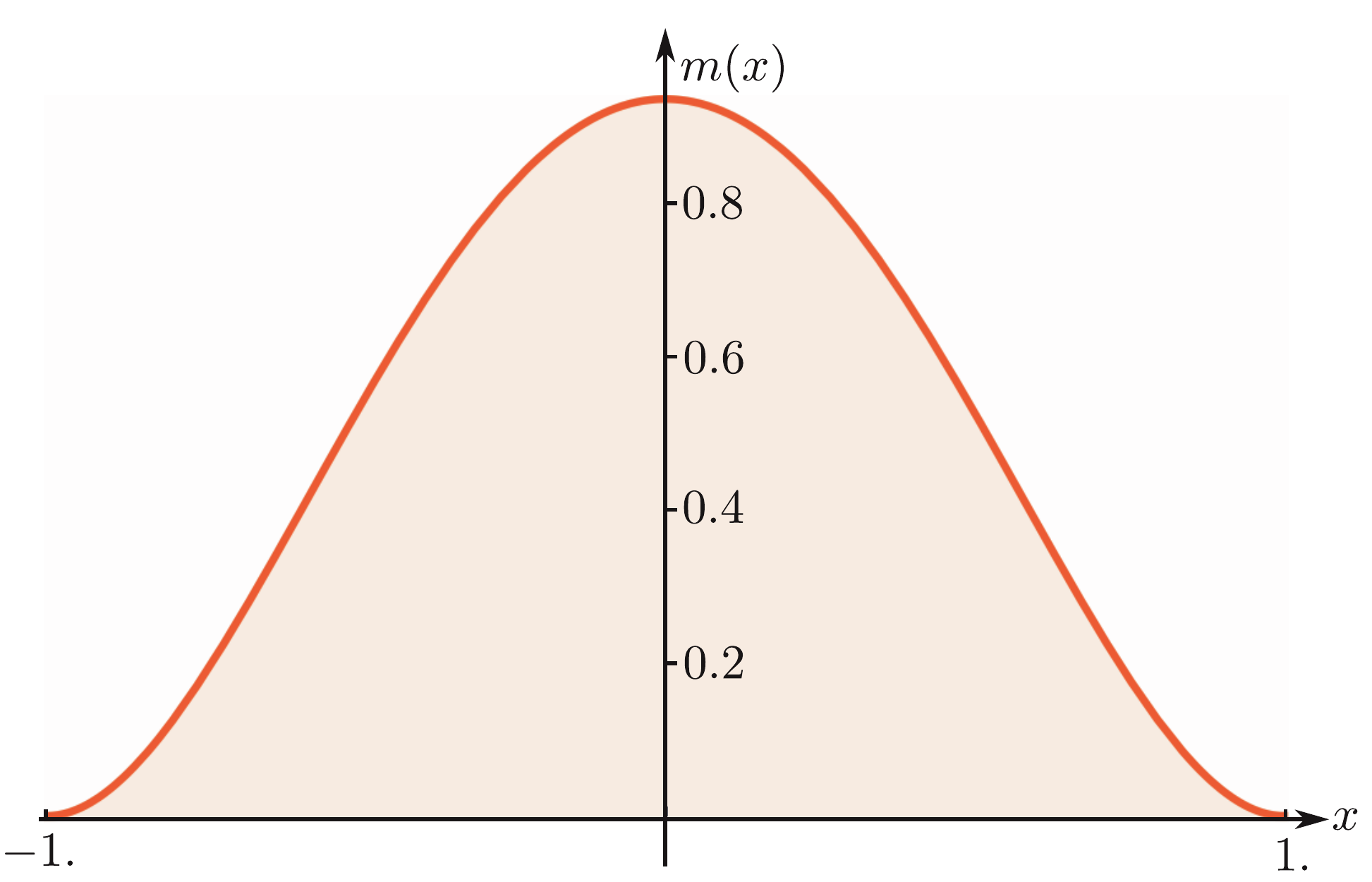}
	\caption{ $C^1$-continuous quartic spline mollifier with the support size $h_m = 2.$}
	\label{fig:kern1D}
\end{figure}

%--------------------------------------------------------------------------------
%
Without loss of generality we focus in the following on trivariate basis functions. As in the univariate case the basis functions for a cell~$\Omega_i$ are given by 
\begin{equation} \label{eq:basis2d}
	\vec N_i (\vec x) = \int_{\Omega_i}  m(\vec x- \vec y) \vec p_i  (\vec y) \D \vec y  \, .
\end{equation}
The vector~$\vec p_i(\vec x )$ contains the scaled and shifted multivariate monomial basis of degree~$q^p$.  In this paper, the mollifier 
\begin{equation}
	m(\vec x) = m \left(x^{(1)} \right)  \cdot m \left(x^{(2)} \right) \cdot m \left(x^{(3)} \right)
\end{equation}
is composed of $C^1$-continuous quartic splines 
\begin{equation} \label{eq:quartSpline}
	m(x) =  
	\begin{cases}
		\frac{15}{8h_m} \left(1- 8  \left ( \frac{x}{h_m} \right ) ^2+ 16 \left ( \frac{x}{h_m} \right ) ^4 \right)  \quad  & \text{if } \, | x |  <   h_m/2 \\
		0						 &  \text{if }  \, | x |  \ge  -  h_m/2
	\end{cases} \, ;
\end{equation}
see Figure~\ref{fig:kern1D}. It is easy to verify that~$m(-h_m/2) = m(h_m/2) = 0$ and $\D m(-h_m/2)/\D x = \D m(h_m/2) / \D x = 0$. Hence, the mollifier is $C^1$ continuous. The continuity of the mollifier can be increased by forcing more derivatives to be zero at~$x=- h_m/2$ and~$x=h_m/2$, which can be achieved either by choosing a higher order polynomial, using a non-polynomial function or introducing more breakpoints.  Obviously all of these approaches increase the cost of evaluating the convolution integrals.  

\begin{figure}[tb]
\centering
	\subfloat[][Cell~$\Omega_i$ and mollifier~$\Box_{\vec x}$] {
		\includegraphics[width=0.28 \textwidth]{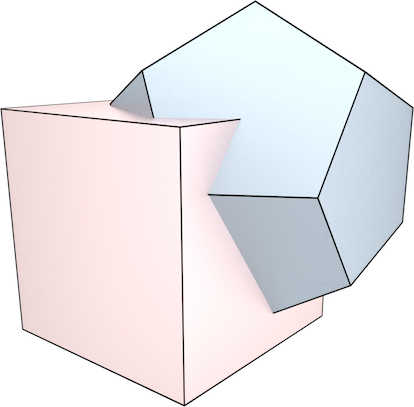} } 
		\hspace{0.05\textwidth}
	\subfloat[][Triangulated  $\omega_i$  \label{fig:basis2Db} ] {
		\includegraphics[width=0.28 \textwidth]{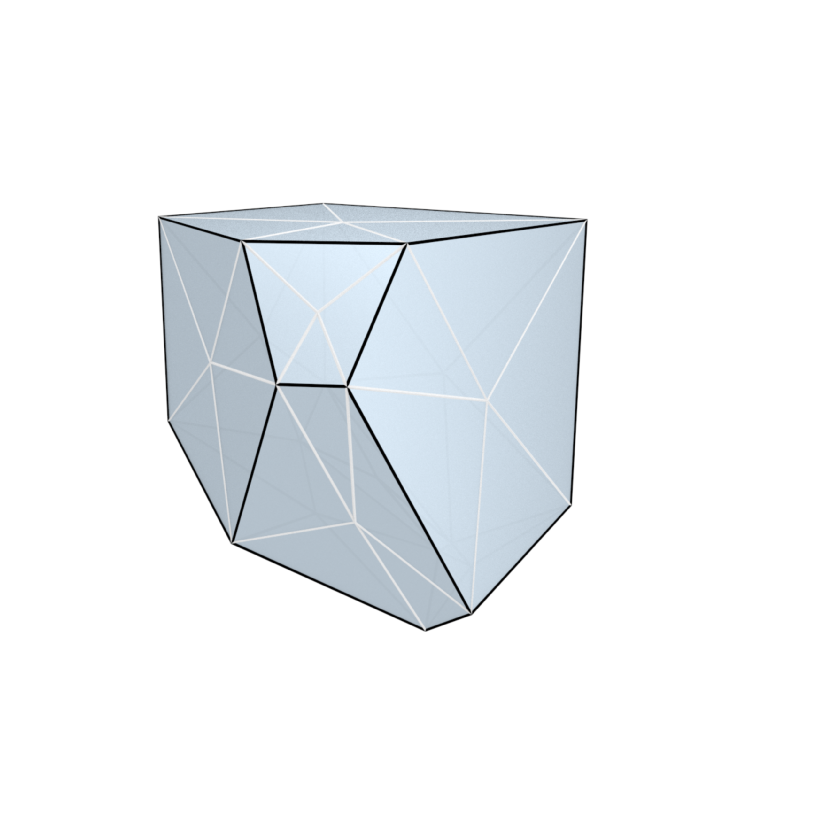} }
	\caption{Evaluation of mollified basis functions~$\vec N_i (\vec x)$ by numerically computing the convolution integral~\eqref{eq:basis2d} at a point~$\vec x \in \Omega$. (a) The support of the mollifier~$\Box_{\vec x}$ centred at~$\vec x$ (in red) and the cell $\Omega_i$ (in blue). (b) The integration domain for the convolution integral is the boolean difference $\omega_i  = \Box_{\vec x} \cap  \Omega_i$. \label{fig:basis2D}}
\end{figure}

The availability of efficient evaluation techniques for the multi-dimensional convolution integral~\eqref{eq:basis2d} is vital to the proposed approach. We note that the integrand is a polynomial and that it is only non-zero on the intersection, i.e. boolean intersection, of the support of the mollifier and the considered cell, i.e., 
\begin{equation}
	\omega_i  \coloneqq \Box_{\vec x} \cap  \Omega_i  \, ,
\end{equation}
where~$\Box_{\vec x} =  \supp m(\vec x - \vec y)$ denotes the support of the mollifier centred at the evaluation point~$\vec x$, see Figure~\ref{fig:basis2D}.  The intersection domain~$\omega_i$  is convex because both ~$\Box_{\vec x} $  and~$\Omega_i$ are convex. Computing the intersection of polytopes and the integration of polynomials over polytopes are recurring tasks in computer graphics and many robust and efficient algorithms and implementations are available. In~\ref{app:intersection} we introduce one such algorithm for determining the intersection between a cell and a mollifier. One possible approach to evaluate the integrals on~$\omega_i$ is first to tesselate it and then to integrate over the obtained simplices using Gaussian quadrature. Considering that the integrands are polynomials and the tessellation consists of affinely mapped simplices, the integration can be performed exactly (up to round-off error) using a sufficient number of quadrature points. 
A more elegant approach is to reduce the domain integrals to line integrals by repeated application of the divergence theorem, see e.g.~\cite{mirtich1996fast}. We briefly sketch the conversion of volume integrals to surface integrals for completeness.  The polytope~$\omega_i$ consists of a set of uniquely orientated faces~$ \{ \gamma_{i,j} \}$, i.e. all the respective normals~$\vec n_{i,j}$ point outside the domain.  To integrate an arbitrary polynomial~$f(\vec x)$ it is first integrated, e.g. in the~$x^{(1)}$ direction, 
\begin{equation}
	\tilde f(\vec x) = \int f(\vec x) \, \D x^{(1)} \, .
\end{equation} 
The divergence theorem applied to this new function yields
\begin{equation} \label{eq:intgFace}
	f(\vec x) = 
	\int_{\omega_i}  \nabla \cdot  
	\begin{pmatrix}
	\tilde f(\vec x) \\
	0 \\
	0
	\end{pmatrix} \D \omega_i
	= \sum_j   n^{(1)}_{i, j} \int_{\gamma_{i, j}} \tilde f (\vec x)   \D \gamma_{i,j} \, ,
\end{equation}
where the surface normal~$\vec n_{i,j}$ of the face~$\gamma_{i,j}$ is constant. It is possible to stop at this point and to numerically evaluate the surface integrals after triangulating the faces~$\gamma_{i,j}$. The number of quadrature points on each face is chosen to integrate exactly polynomials of degree $q^m + q^p$. However, it is also possible to reapply the divergence theorem to reduce the surface integrals to line integrals, which can then be analytically evaluated. In our implementation we evaluate~\eqref{eq:intgFace} numerically by triangulating the faces $\gamma_{i,j}$, as indicated in Figure~\ref{fig:basis2Db}. Note that the change from volume to surface integrals already yields a significant reduction in the number of integration points. 

\begin{figure}
	\centering
	\includegraphics[width=0.6 \textwidth]{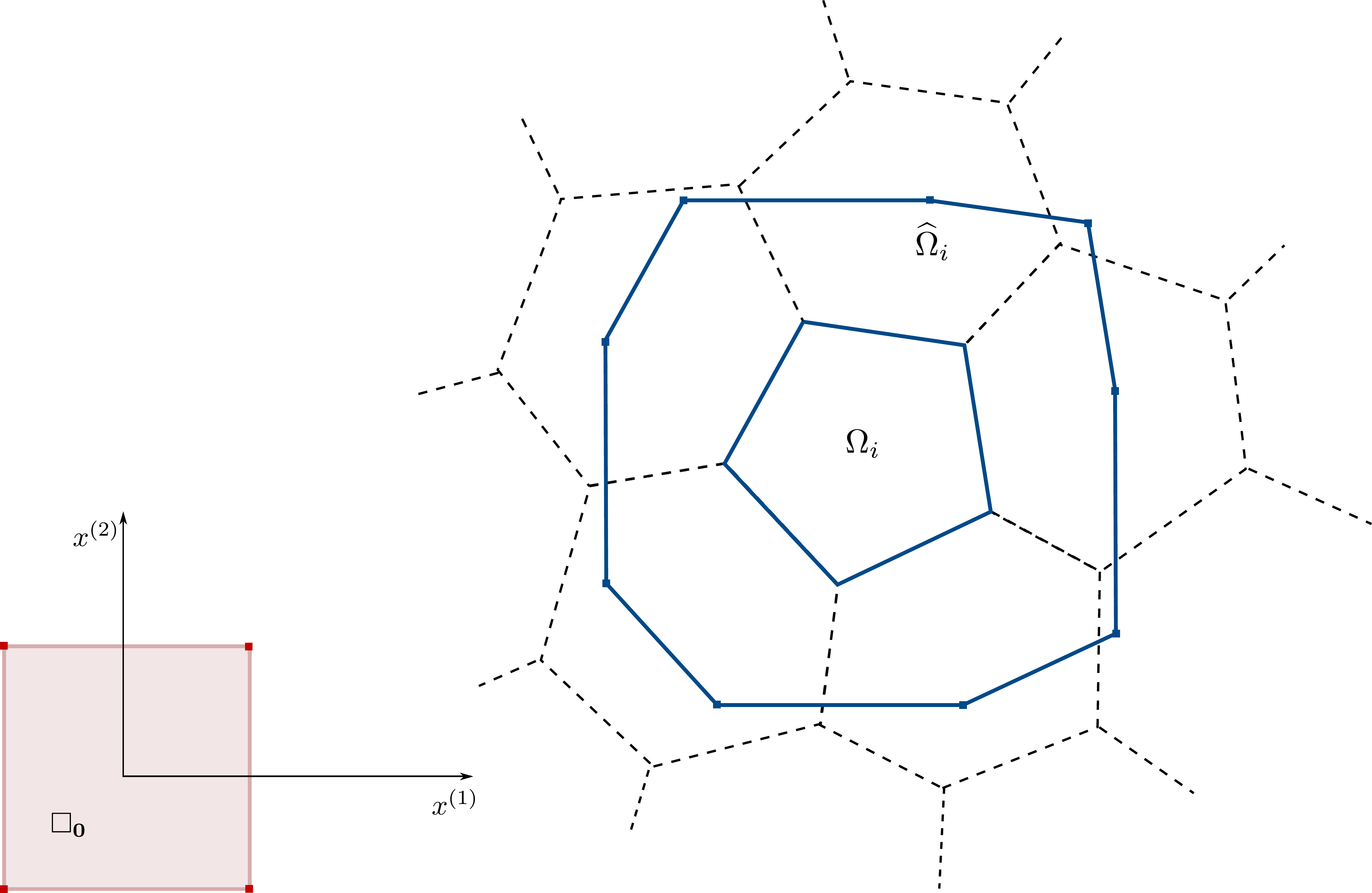} 
	\caption{A cutout from a Voronoi tesselation and the support of the mollifier. The mollified basis functions~$\vec N_i (\vec x)$ of the highlighted cell~$\Omega_i$ are supported on the larger polygonal domain~$\widehat \Omega_i = \supp \vec N_i (\vec x)$.  The polygonal domain is the Minkowski sum, i.e.~$\widehat \Omega_i = \Omega_i \oplus \Box_{\vec 0}$,  of the cell~$\Omega_i$ with the mollifier support~$\Box_{\vec 0}$. \label{fig:basisSupp2D}}
\end{figure}

As in the univariate case, the support of a mollified basis function~$\vec N_i(\vec x)$ is larger than its respective cell, see Figure~\ref{fig:basisSupp2D}. In finite element computations also the support of a basis function is required, which is given by the Minkowski sum   
\begin{equation}
	\widehat \Omega_i = \Omega_i \oplus \Box_{\vec 0} =  \{ \vec y + \vec z | \,   \vec  y \in \Box_{\vec 0}, \, \vec z \in \Omega_i  \}   \, .
\end{equation}
In~\ref{app:minkowski} we introduce the algorithm used for computing the Minkowski sum of two convex polyhedra.

%
%--------------------------------------------------------------------------------          
\section{Finite element discretisation with mollified basis functions \label{sec:fem}}
%--------------------------------------------------------------------------------
%
The smoothness and approximation properties of the mollified basis functions make them ideal for finite element analysis.  In the following, we briefly outline the discretisation of a Poisson equation using the mollified basis functions. As in the previous sections we assume that a partitioning of the domain consisting of convex polytopes is given. We generate such a partitioning using a Voronoi diagram from a given implicit, i.e. level set, or parametric, i.e. surface mesh, description of the domain boundary. A set of points is placed within and outside the domain to generate the Voronoi diagram. The points outside the domain ensure that the mollifier's support is fully covered by local polynomials when the mollifier is placed on the domain boundary. The Voronoi diagram and the respective mollified basis functions do not conform to the domain boundaries. Therefore, the boundary conditions are applied weakly on cells cut by the boundary, as in immersed or embedded finite elements. 

Notwithstanding this, it appears to be possible to continuously shrink the mollifier support when approaching the boundaries from the inside. In the limit on the boundary the mollifier becomes a Dirac delta and the interpolation within the domain becomes independent of the outside of the domain. This idea has, however, not been further pursued in this paper in order to focus on other aspects of the method. 

The Poisson equation on a domain~$\Omega$ is given by
\begin{equation}\label{eq:femPoissonStrongForm}
\begin{aligned}
   -\nabla \cdot \nabla{u}
&  = s 
&& \text{in $\Omega$}
\\
   u
&  = \overline{u}
&& \text{on $\Gamma_D$}
\\
   \vec{n} \cdot \nabla{u}
&  = \overline{t}
&& \text{on $\Gamma_N$}
\end{aligned}
\end{equation}
where $\overline{u}$ is the prescribed solution field on the Dirichlet boundary~$\Gamma_D$ and~$\overline{t}$ is the prescribed flux on the Neumann boundary~$\Gamma_N$ with the outward normal~$\vec{n}$. 
The  weak formulation of the Poisson equation can be stated according to Nitsche~\cite{Nitsche:1971aa} as: Find $u \in H^1(\Omega)$ such that
\begin{equation}\label{eq:femWeakForm}
   \underbrace{
      \int_\Omega \nabla{u} \cdot \nabla{v}\D\Omega
   }_{\displaystyle a(u,v)}
 = \underbrace{
      \int_\Omega s \, v \D\Omega
      + \int_{\Gamma_N} \overline{t} \, v \D\Gamma
   }_{\displaystyle b(v)}
   - \underbrace{
      \gamma \int_{\Gamma_D} (u-\overline{u}) \, v \D\Gamma
   }_{\displaystyle \gamma p(u,v)}
  + \underbrace{
      \int_{\Gamma_D} \Big(
         (u-\overline{u}) \, \vec{n} \cdot \nabla v
         + (\vec{n} \cdot \nabla u)\, v
      \Big) \D\Gamma
   }_{\displaystyle l(u,v)} \, . 
\end{equation}
The stabilisation parameter can be set to~$\gamma=0$ when the sign of the last term~$(\vec n \cdot \nabla u) v$ is reversed, as proposed in a number of papers~\cite{oden1998discontinuoushpfinite, boiveau2015penalty, schillinger2016non}. In our computations we use this so-called non-symmetric Nitsche method which has a non-symmetrical system matrix. The trial and test functions are discretised with the mollified basis functions 
\begin{equation}\label{eq:femBasis}
   u^h(\vec {x})   = \sum_{i=1}^{n_c} \vec{N}_i (\vec {x}) \cdot  \vec {\alpha}_i  \quad   \text{and} \quad  v^h(\vec{x})  = \sum_{j=1}^{n_c}   \vec{N}_j (\vec {x}) \cdot  \vec {\beta}_j \, ,
\end{equation}
where~$n_c$ is the number of the polytopic cells in the mesh.

Introducing the interpolation equations~\eqref{eq:femBasis} into the weak form~\eqref{eq:femWeakForm} yields a linear system of equations with the unknowns~$\vec \alpha_i$,  which are the coefficients of the local monomial bases in the cells. For instance, the bilinear form~$a(u^h, \, v^h)$ becomes after discretisation 
\begin{equation}
	a(u^h, \, v^h) = \sum_{i=1}^{n_c} \sum_{j=1}^{n_c}  \vec \alpha_i^\trans  \int_\Omega \nabla \vec N_i   \nabla \vec N_j^\trans \D \Omega \, \, \vec \beta_j \, .
\end{equation}
As usual, the domain integral is evaluated numerically  after splitting it into cell contributions 
\begin{equation} \label{eq:weakPartitioned}
	a(u^h, \, v^h) = \sum_{k=1}^{n_c} \left ( \sum_{i=1}^{n_c} \sum_{j=1}^{n_c}  \vec \alpha_i^\trans  \int_{\Omega_k} \nabla \vec N_i \nabla \vec N_j^\trans \D \Omega_k \, \, \vec \beta_j \right ) \, .
\end{equation}
We evaluate the integral over a cell~$\Omega_k$ by first decomposing it into tetrahedra and then applying standard Gauss integration. A cell is tetrahedralised by introducing additional nodes at its centre and face centres.  At each integration point the mollified basis functions are evaluated as described in Section~\ref{sec:basis}. Although the sketched integration of the weak form~\eqref{eq:weakPartitioned} is straightforward, Gauss integration unfortunately requires too many quadrature points because of the breakpoints of the basis functions within the cells. Using too few quadrature points usually leads to suboptimal convergence rates. 

The principal difficulties encountered in efficient and accurate integration of the weak form~\eqref{eq:weakPartitioned} are very similar to those encountered in meshless methods. In the variationally consistent integration techniques for meshless methods the test functions are modified to satisfy a consistency condition when integrated numerically, see e.g.~\cite{hillman2013, hillman2015}. If the basis functions can reproduce a solution~$\vec u (\vec x)$ of polynomial degree~$q^p$, the finite element scheme must be able to solve exactly such a problem even when the integrals are evaluated numerically. Assuming a problem with the solution of polynomial degree~$q^p$ and inserting it into the discretised weak form yields for each cell a consistency condition for integration. To satisfy this consistency condition in a cell~$j$ the gradient of the test functions is modified as 
\begin{equation}
	\nabla  \widetilde{\vec {  N}}_j^\trans  (\vec x) = \nabla \vec N_j^\trans  (\vec x)+ \vec \Lambda_j  \vec \psi_j (\vec x) \otimes \vec 1 \, ,
\end{equation}
where~$\vec 1$ is a vector of all ones and the vector~$ \vec \psi_j (\vec x)$ contains a monomial basis of degree~$q^{p}-1$ with the yet to be determined matrix of coefficients $\vec \Lambda_j$. The support of the basis functions~$\vec \psi_j (\vec x)$ is chosen to be same as of the basis functions~$ \vec N_j (\vec x)$.  The matrix of coefficients~$\vec \Lambda_j$ is determined by solving per basis function one small linear equation system so that the variational consistency condition is satisfied. This equation system contains the integrals of the basis functions over their supports, which is given by the Minkowski sum of the mollifier support and the cell. See~\cite{hillman2013} for further details. 

In the cells cut by the domain boundary the element integrals are evaluated only over parts of the cell which lie inside the domain. The respective integration domains are obtained by clipping the cells, see \ref{app:clipVoro}. The resulting polyhedron is tetrahedralised with the same approach used for a non-clipped cell. As in immersed, or embedded domain, methods the faces of the tetrahedra can be projected to the curved domain boundaries if higher-order boundary approximation is desired, see e.g.~\cite{cheng2010higher, kudela2015efficient, xiao2019immersed}.

%
%--------------------------------------------------------------------------------          
\section{Examples \label{sec:examples}}
%--------------------------------------------------------------------------------
%
We introduce in this section several examples of increasing complexity to experimentally verify the convergence of the proposed mollified finite element approach. In the one-dimensional problems both the convolution and finite element integrals are evaluated analytically, whereas in multi-dimensional problems both integrals are evaluated numerically.  In addition, in multi-dimensional problems only the quartic spline mollifier~\eqref{eq:quartSpline} consisting of a single polynomial with no internal breakpoints is used.  
The mollified basis functions contain monomials of up to degree $q^p + q^m + 1$, although they are only complete up to degree $q^p$, and are non-zero over several cells. Therefore, it is not obvious how many quadrature points to choose in each of the integration triangles used for evaluating the finite element integrals. In our present computations we determine a stable number of quadrature points by successively increasing their number until we have a stable solution. In two dimensional problems we choose for linear basis functions ($q^p =1$) three integration points for domain integrals and five for boundary integrals, and for quadratic ($q^p=2$) we choose four and five integration points respectively. In all problems the Dirichlet boundary conditions are enforced with the parameter-free non-symmetric Nitsche method.

%
%--------------------------------------------------------------------------------          
\subsection{One-dimensional Poisson problem}
%--------------------------------------------------------------------------------
%
As a first example we consider the  solution of the one-dimensional Poisson-Dirichlet problem $- \D^2 u / \D x^2  =  s$ on the domain $\Omega = (0, \, 1)$. The source term~$s(x)$ is chosen such that the solution is equal to 
\begin{equation}
	 u(x) = \sin (3 \pi x).
\end{equation} 
The initial coarse mesh consisting of~$n_c=6$ cells is chosen to be non-uniform. The cell sizes, starting from the left, are  $h_{c,1} = 0.15$, $h_{c,2} = 0.2$, $h_{c,3} = 0.15$, $h_{c,4} = 0.15$, $h_{c,5} = 0.2$  and $h_{c,6} = 0.15$. In addition to these cells, each domain boundary is padded with an extra ghost cell to ensure that the obtained mollified basis functions have the same approximation properties over the entire domain. We obtain  finer meshes by repeated bisectioning of all cells.

\begin{figure} 
	\centering
	\subfloat[$L^2$-norm error]{\includegraphics[width=0.52\textwidth]{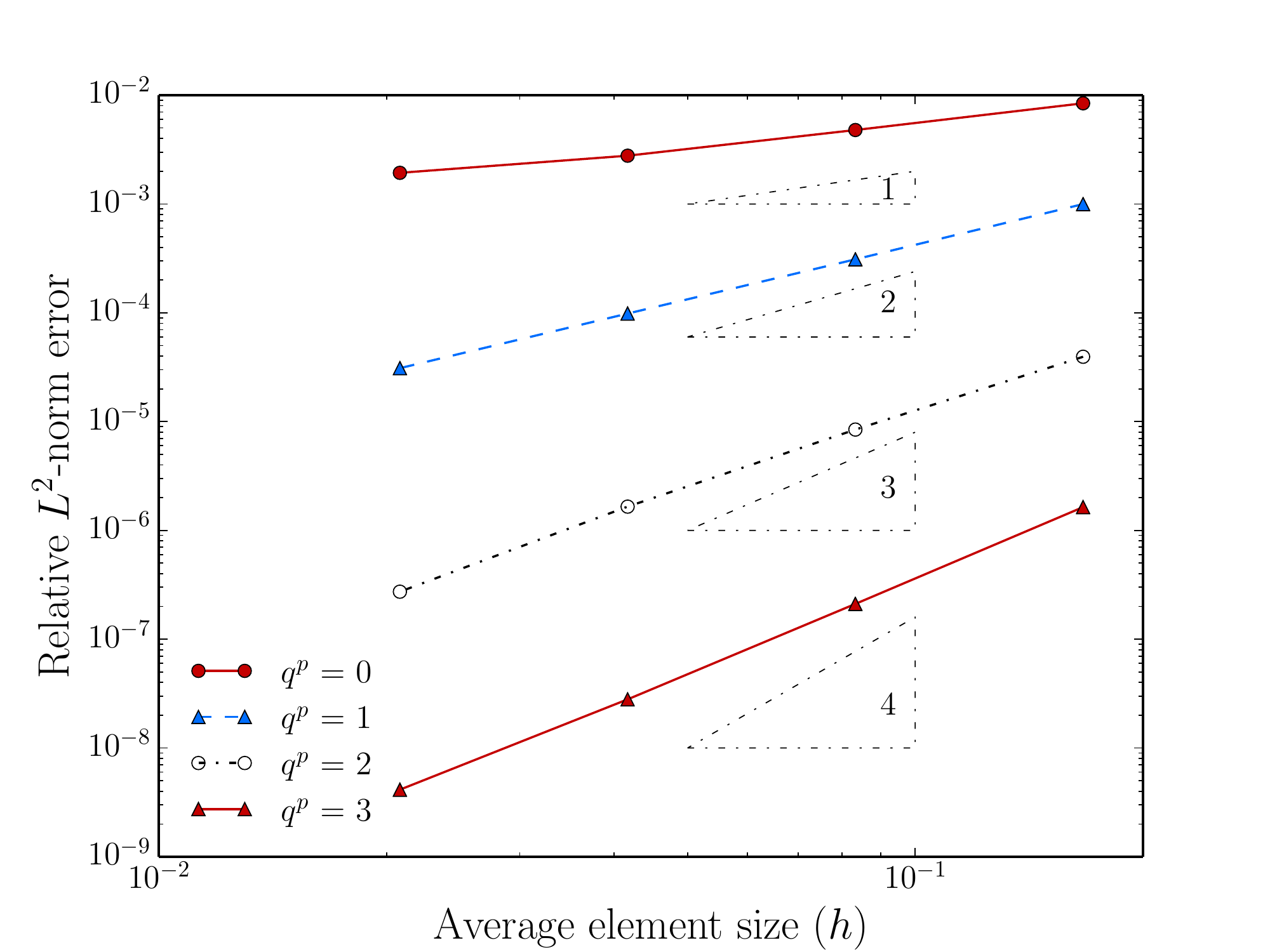}}
	\subfloat[$H^1$-seminorm error]{\includegraphics[width=0.52\textwidth]{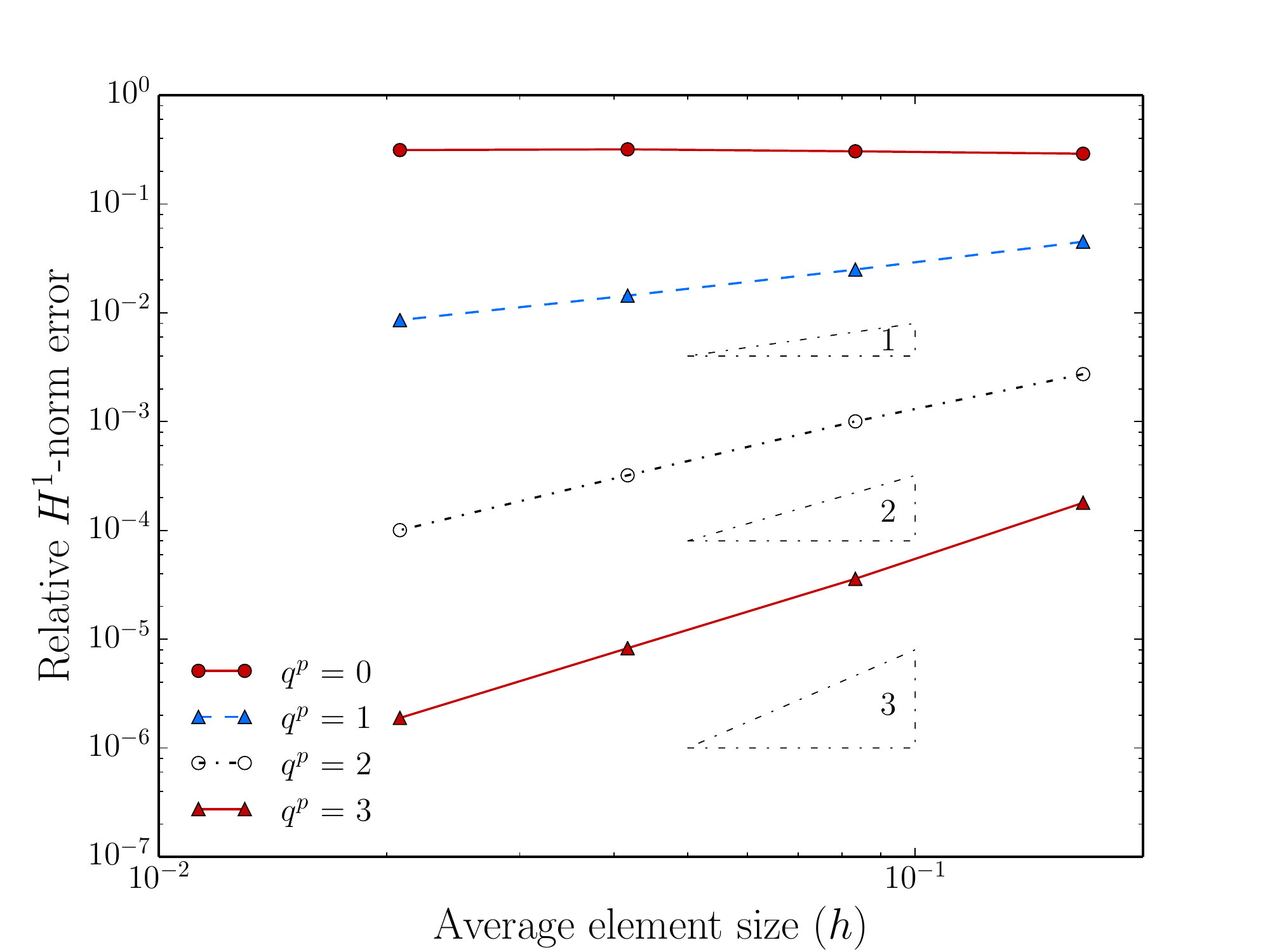}}
	\caption{One-dimensional Poisson problem. Convergence with normalised linear B-spline mollifier and local polynomial basis of degrees \mbox{$q^p \in \{ 0, 1,2,3\}$.} }
	\label{fig:conv1DvarP}
\end{figure}

In the following set of experiments we study the influence of the choice of the local polynomial basis $\vec p_i(x)$ and the mollifier~$m(x)$ on the convergence of the finite element solution. Firstly, we take in turn different local polynomials of degrees~$q^p \in \{ 0, 1,2,3\}$ and a normalised linear B-spline, i.e. hat function, mollifier with a support width of 
\begin{equation}
	h_m =  2  \left ( \max_j h_{c,j} \right )  \, .
\end{equation}
Note that the normalisation of the mollifier is essential for ensuring that it integrates as required to one.  Figure~\ref{fig:conv1DvarP} shows the optimal convergence of the mollified finite element approach in the $L^2$ norm and $H^1$ seminorm with $q^p+1$ and $q^p$ respectively. These convergence rates are in agreement with the analytic estimates provided in~\ref{app:convergence}. 
%This can be explained by the ability of the mollified basis functions exactly to represent polynomials of degree~$q^p$.
%   

\begin{figure}
	\centering
	\subfloat[$L^2$-norm error]{\includegraphics[width=0.5\textwidth]{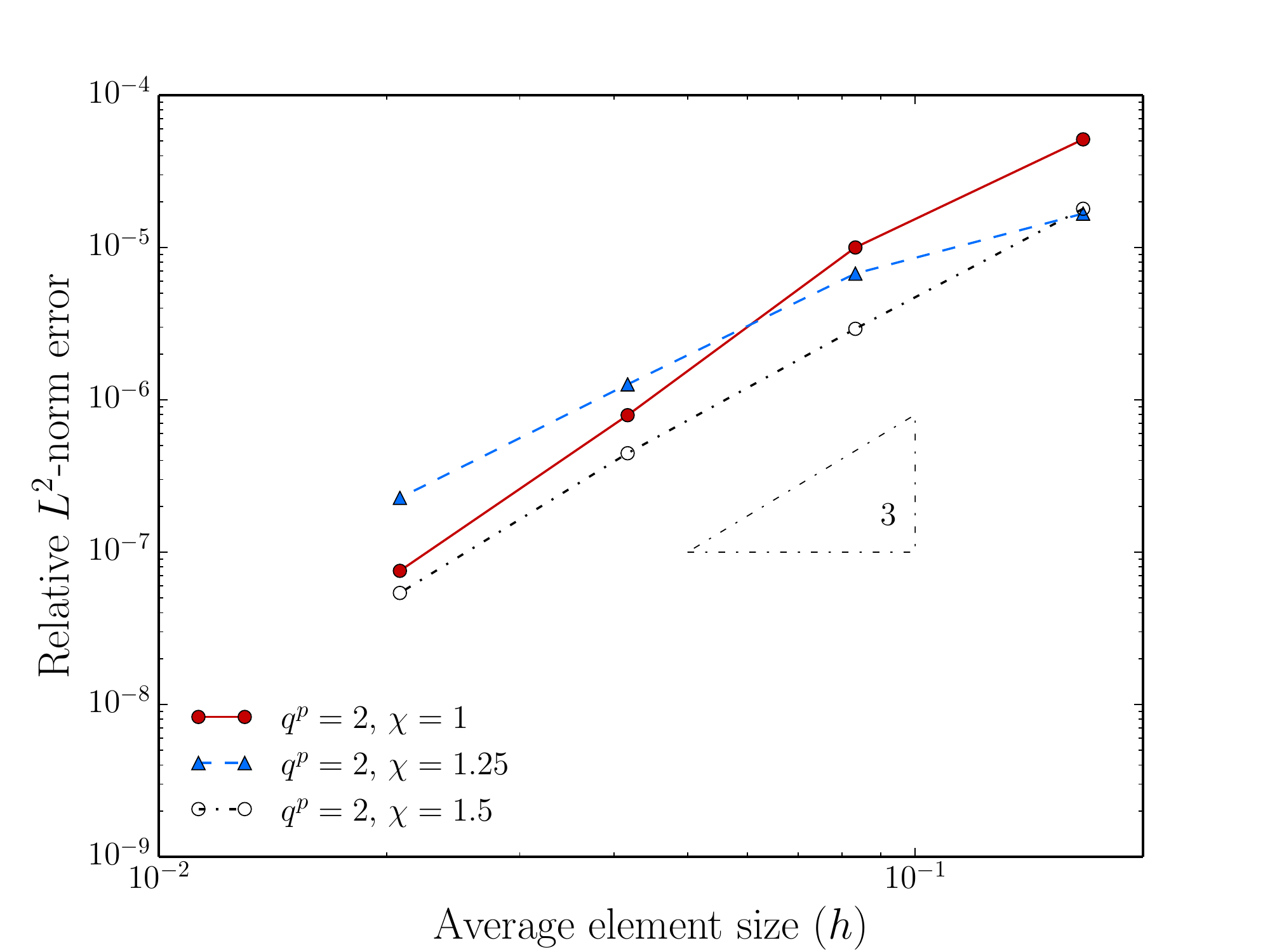}}
	\subfloat[$H^1$-seminorm error]{\includegraphics[width=0.5\textwidth]{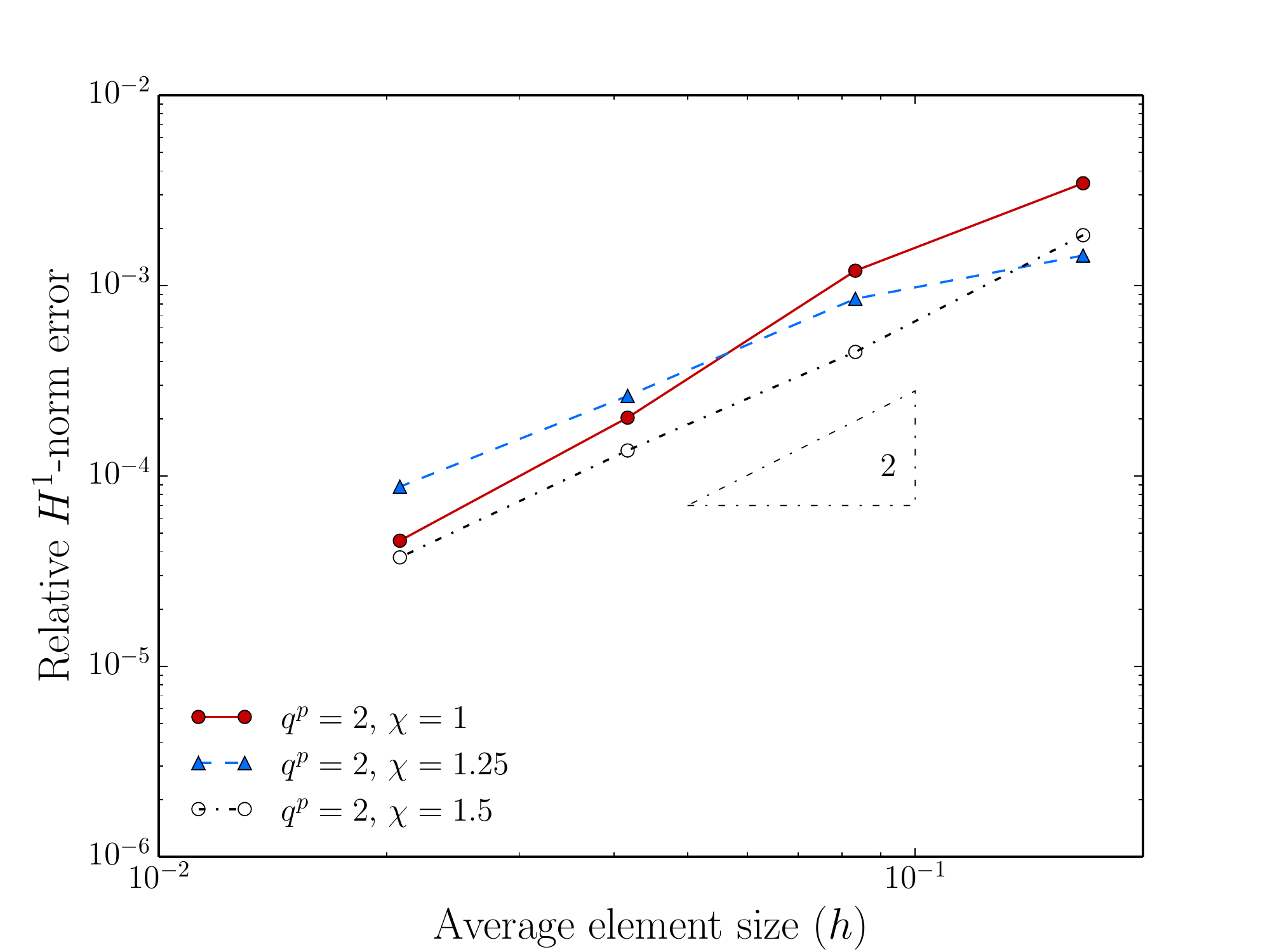}}
	\caption{One-dimensional Poisson problem. Convergence with local polynomial basis of degree $q^p=2$ and normalised linear B-spline mollifier with different support sizes of \mbox{$\chi \in \{1, \,1.25, \,1.5 \}$.} }
	\label{fig:conv1DvarKernSizeQuad}
\end{figure}

Next, we investigate the influence of the mollifier support width on the convergence order while keeping the normalised linear B-spline mollifier. The mollifier width is chosen according to 
\begin{equation}
	h_m =  2  \chi \left ( \max_j h_{c,j} \right )  \quad \text{with } \chi \in \{1, \,1.25, \,1.5 \} \, .  
\end{equation} 
The increase in mollifier size leads to an increase in the support size of the mollified basis functions, which results in an increase of the number of non-zero basis functions in a cell. The obtained optimal convergence rates for~$q^p=2$ are shown in Figure~\ref{fig:conv1DvarKernSizeQuad}. The increase in mollifier width leads to a somewhat decrease in the convergence constants, but the optimal support size appears to depend on the specific problem considered.  The results for higher order local polynomials are similar and have been not included here. 

\begin{figure} 
	\centering
	\subfloat[$L^2$-norm error]{\includegraphics[width=0.5\textwidth]{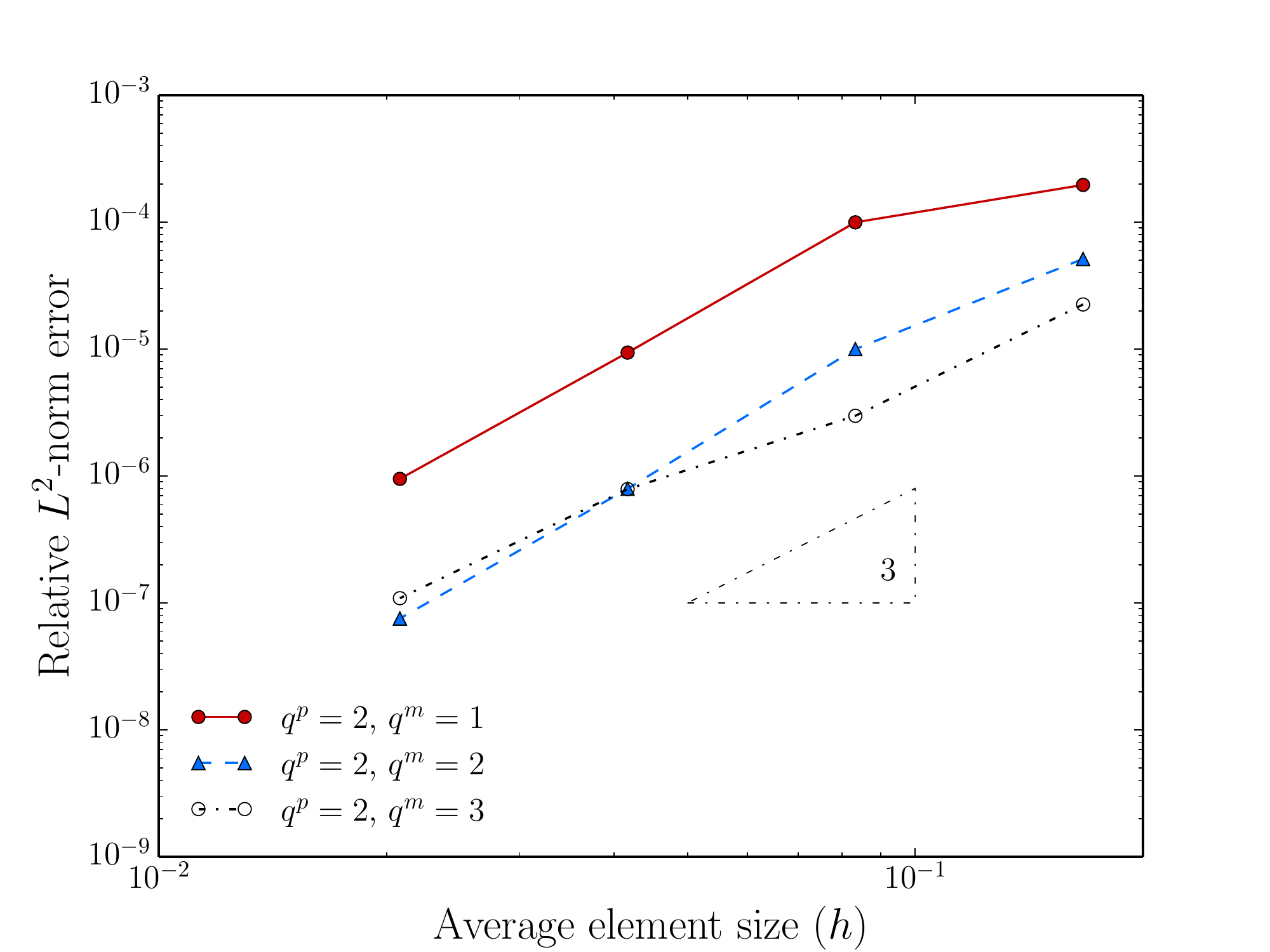}}
	\subfloat[$H^1$-seminorm error]{\includegraphics[width=0.5\textwidth]{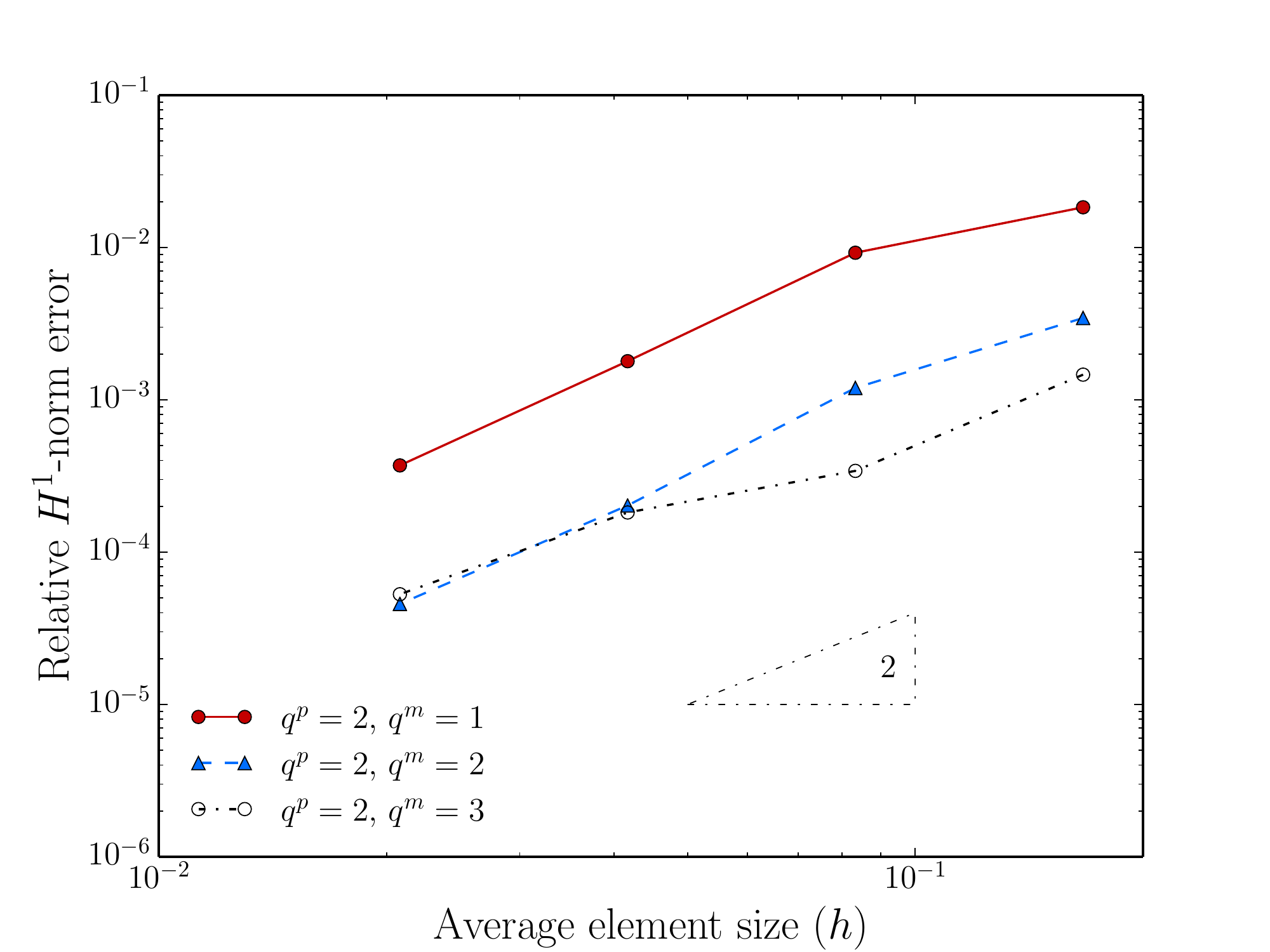}}
	\caption{One-dimensional Poisson problem.  Convergence with local polynomial basis of degree $q^p=2$  and normalised B-spline mollifiers of degrees \mbox{$q^m \in \{ 1, \,2, \,3 \}$.} }
	\label{fig:conv1DvarKernShape}
\end{figure}

Finally, we study the effect of the mollifier smoothness on finite element convergence.  The normalised B-spline mollifiers are of degree $q^m \in \{1,2,3\}$ and the local polynomial is of degree $q^p=2$. The mollifier width factor is chosen as $\chi = 1$. Note that the B-spline mollifiers are~$C^{q^{m}-1}$ continuous so that the obtained mollified basis functions are~$C^{q^{m}}$ continuous. Figure \ref{fig:conv1DvarKernShape} shows that an increase of the mollifier degree and smoothness does not have an effect on the optimal convergence rate, but leads to a significant decrease in convergence constants. 

%
%--------------------------------------------------------------------------------          
\subsection{Two-dimensional examples}
%--------------------------------------------------------------------------------
%

%
%--------------------------------------------------------------------------------          
\subsubsection{Poisson problem on a square domain}
\begin{figure}
	\centering
	\subfloat[$n_c = 16$]{\includegraphics[width=0.3\textwidth]{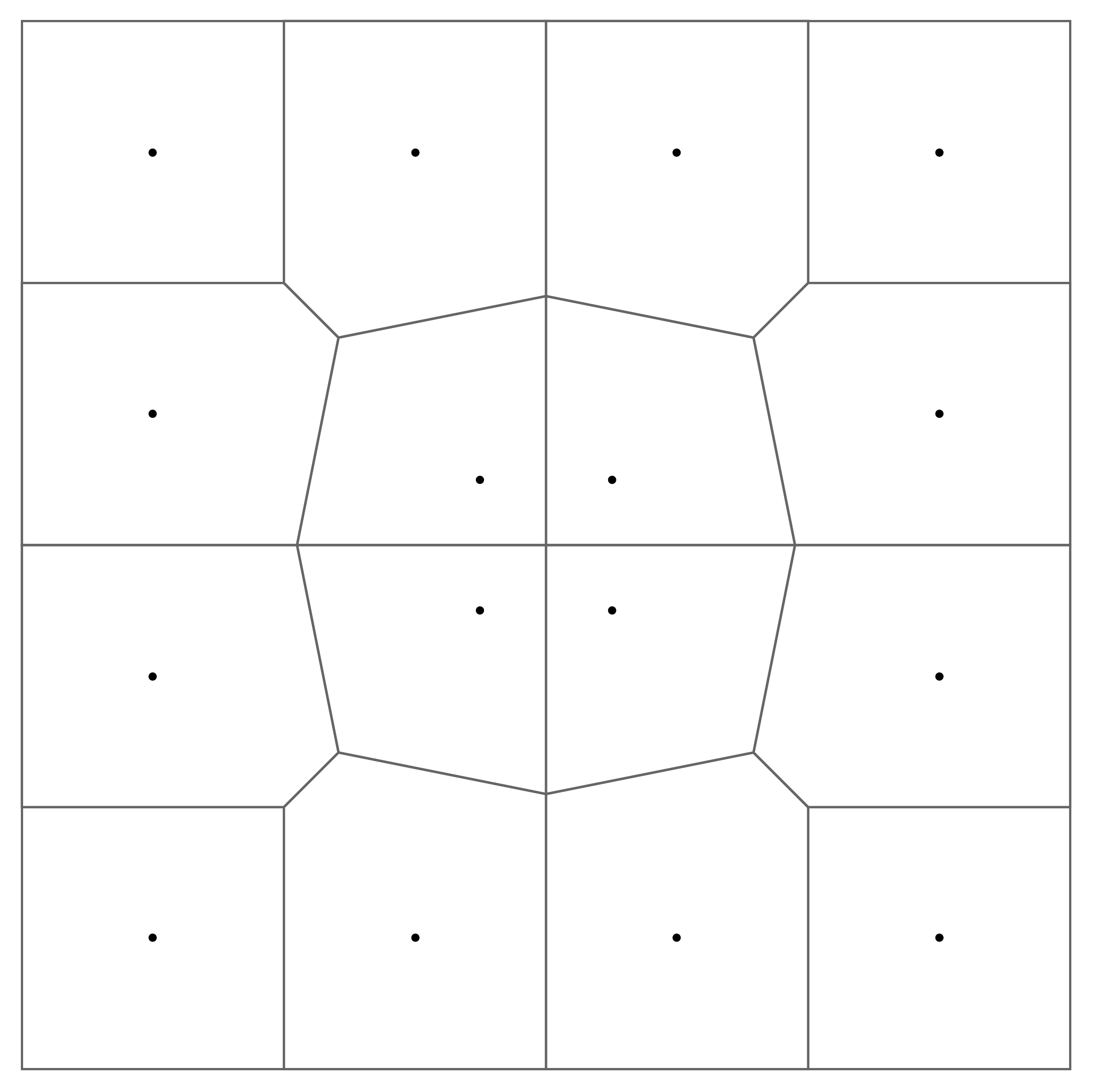}\label{fig:cells2D16}} \qquad
	\subfloat[$n_c = 64$]{\includegraphics[width=0.3\textwidth]{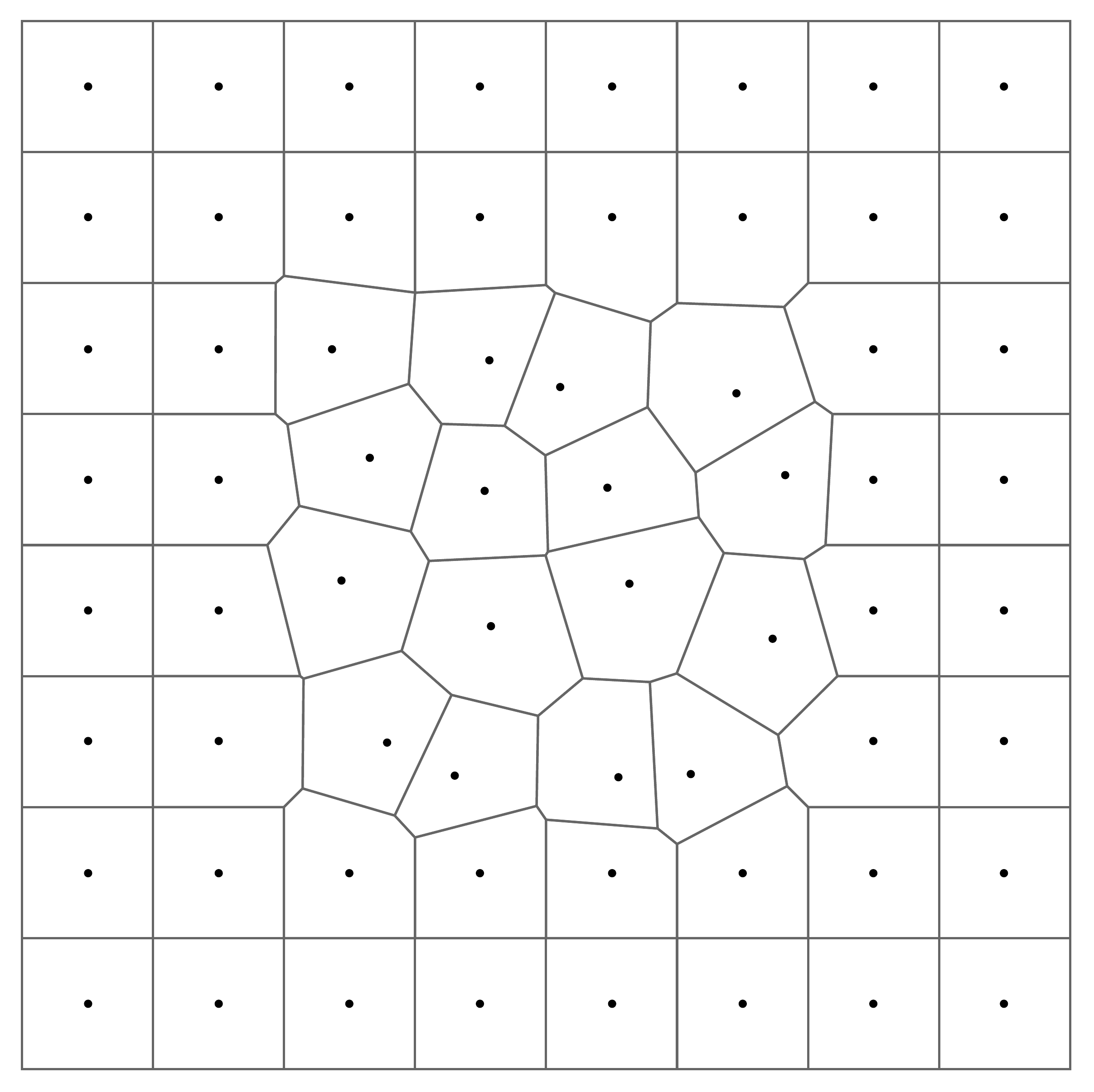}\label{fig:cells2D64}}\qquad
	\subfloat[$n_c = 256$]{\includegraphics[width=0.3\textwidth]{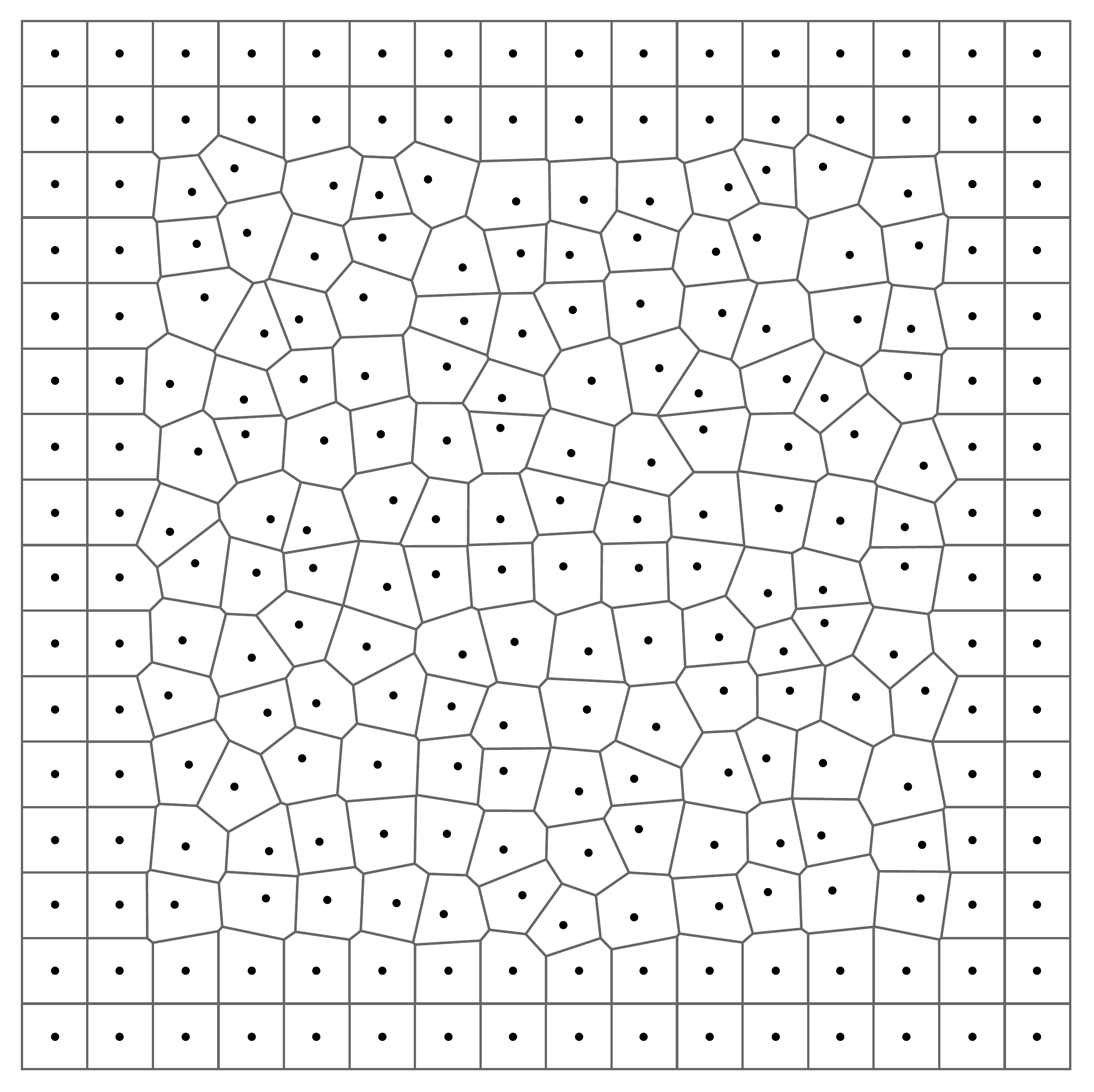}\label{fig:cells2D256}}
	\caption{Poisson problem on a square domain.  Three different partitionings of the domain using the Voronoi diagram of~$n_c$  non-uniformly distributed points.}
	\label{fig:cells2D}
\end{figure}

%--------------------------------------------------------------------------------
%
We consider next the Poisson-Dirichlet problem $-\nabla \cdot \nabla u = s$ on a square domain~$\Omega = ( 0 , \, 1 ) \times ( 0 , \, 1 )$.  The domain~$\Omega$ is partitioned into $n_c$ cells using the Voronoi diagram of~$n_c$  non-uniformly distributed points, see Figure~\ref{fig:cells2D}. Starting from a set of uniformly distributed points we introduce non-uniformity by randomly perturbing their coordinates by $\kappa \sqrt{1/n_c}$ with~$\kappa \sim \set U ( - 0.15 h_m/2, + 0.15 h_m/2)$. Only the coordinates of points farther than a certain distance from the boundaries are perturbed. 

As in the one-dimensional case the domain is padded with an extra layer of ghost cells (not shown in Figure~\ref{fig:cells2D}) to ensure that the mollified basis functions are complete close to the boundaries. Depending on the number of Voronoi cells~$n_c$ the width of the $C^1$ continuous quartic spline mollifier is chosen with 
\begin{equation}
	h_m =  2 \, \bigg( \dfrac{1}{n_c} \bigg)^{\frac{1}{2}} \, .
	\label{eq:sizekern2D}
\end{equation}

We firstly perform a patch test to verify that the mollified finite element method in combination with variationally consistent integration can exactly solve problems with polynomial degree~$q^p$. To this end, we consider on the mesh shown in Figure~\ref{fig:cells2D64} two problems with the exact solutions~\mbox{$u (\vec x) = x + 2 y$} and \mbox{$\ u(\vec x) = x + 2 y + x^2 + 2 x y + y^2$.} Solving the linear problem using the linear mollified basis functions with~$q^{p}=1$ leads to an $L^2$ norm  error of \mbox{$ 8.635 \times 10^{-14} $} and an $H^1$ seminorm error of \mbox{$ 2.902 \times 10^{-12} $.} The corresponding errors for the quadratic problem using quadratic mollified basis functions with~$q^{p}=2$  are \mbox{$ 7.060 \times 10^{-12} $} and \mbox{$ 3.198 \times 10^{-10} $.} This clearly confirms that the mollified finite element method satisfies the patch test. 
 
\begin{figure} 
	\centering
	\subfloat[$L^2$-norm error]{\includegraphics[width=0.5\textwidth]{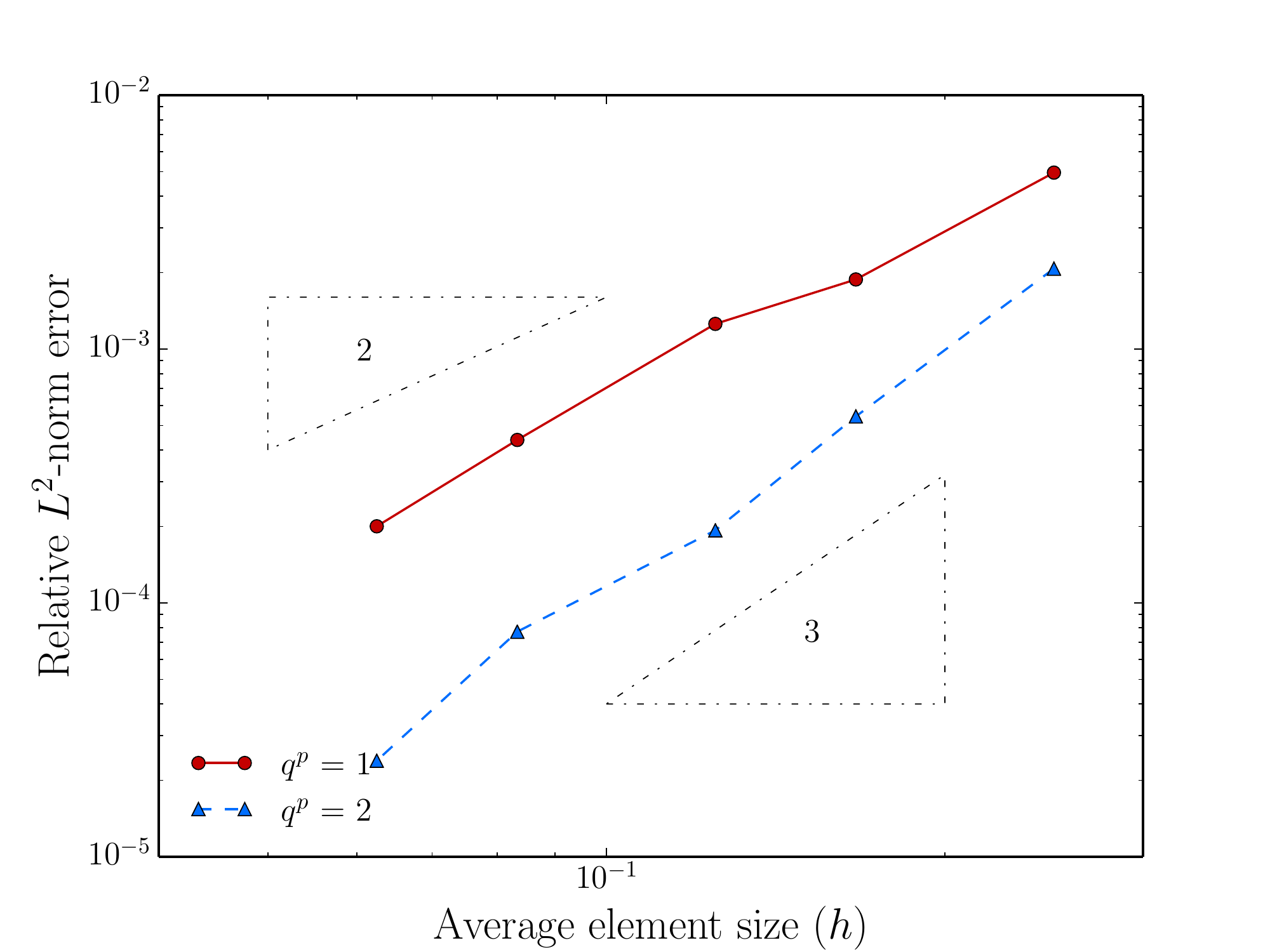}}
	\subfloat[$H^1$-seminorm error]{\includegraphics[width=0.5\textwidth]{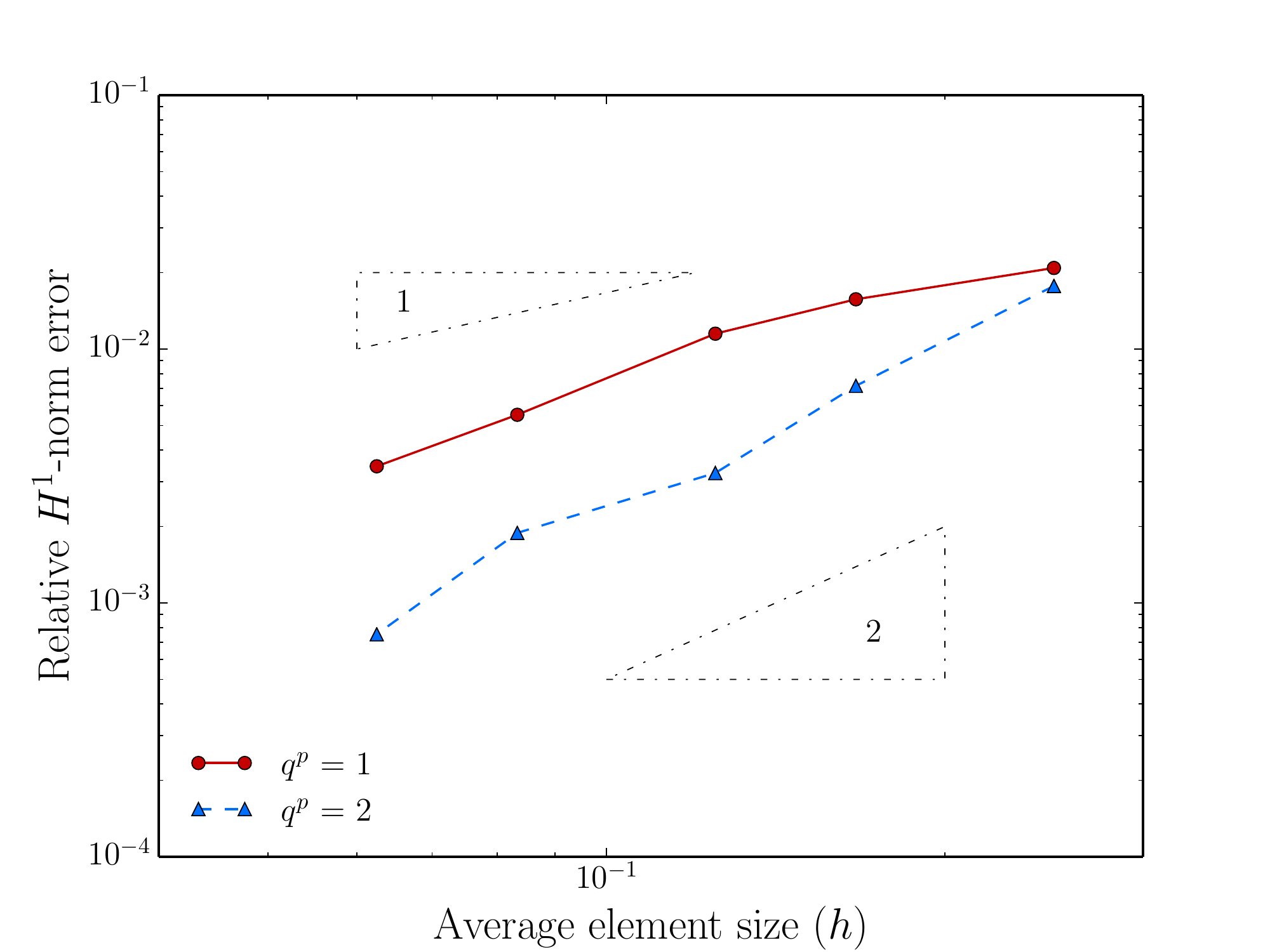}}
	\caption{Poisson problem on a square domain. Convergence with quartic spline mollifier and local polynomial basis of degree $q^{p}=1$ and $q^{p}=2$.  }
	\label{fig:conv2DPoi}
\end{figure}

With the consistency of the method confirmed, we proceed to establish its convergence under mesh refinement. The source term~$s(x)$ is now chosen such that the solution is equal to
 \begin{equation}
 	u(\vec x) = \sin \left (\pi x^{(1)} \right )   \sin \left ( \pi x^{(2)} \right )  \, .
 \end{equation}
The used mollified basis functions are the $C^2$-continuous linear and quadratic basis functions with~$q^{p}=1$ and~$q^{p}=2$, respectively. Figure~\ref{fig:conv2DPoi} shows the convergence of the errors in $L^2$ norm and $H^1$ seminorm as the mesh is refined. Note that the refined meshes are not nested so that some small kinks in the convergence curves may be expected. The average convergence rates are, however, close to optimal, as indicated by the dashed triangles in Figure \ref{fig:conv2DPoi}. 

%
%--------------------------------------------------------------------------------          
\subsubsection{Elastic plate with a hole}
\begin{figure} 
	\centering
	\includegraphics[width=0.375\textwidth]{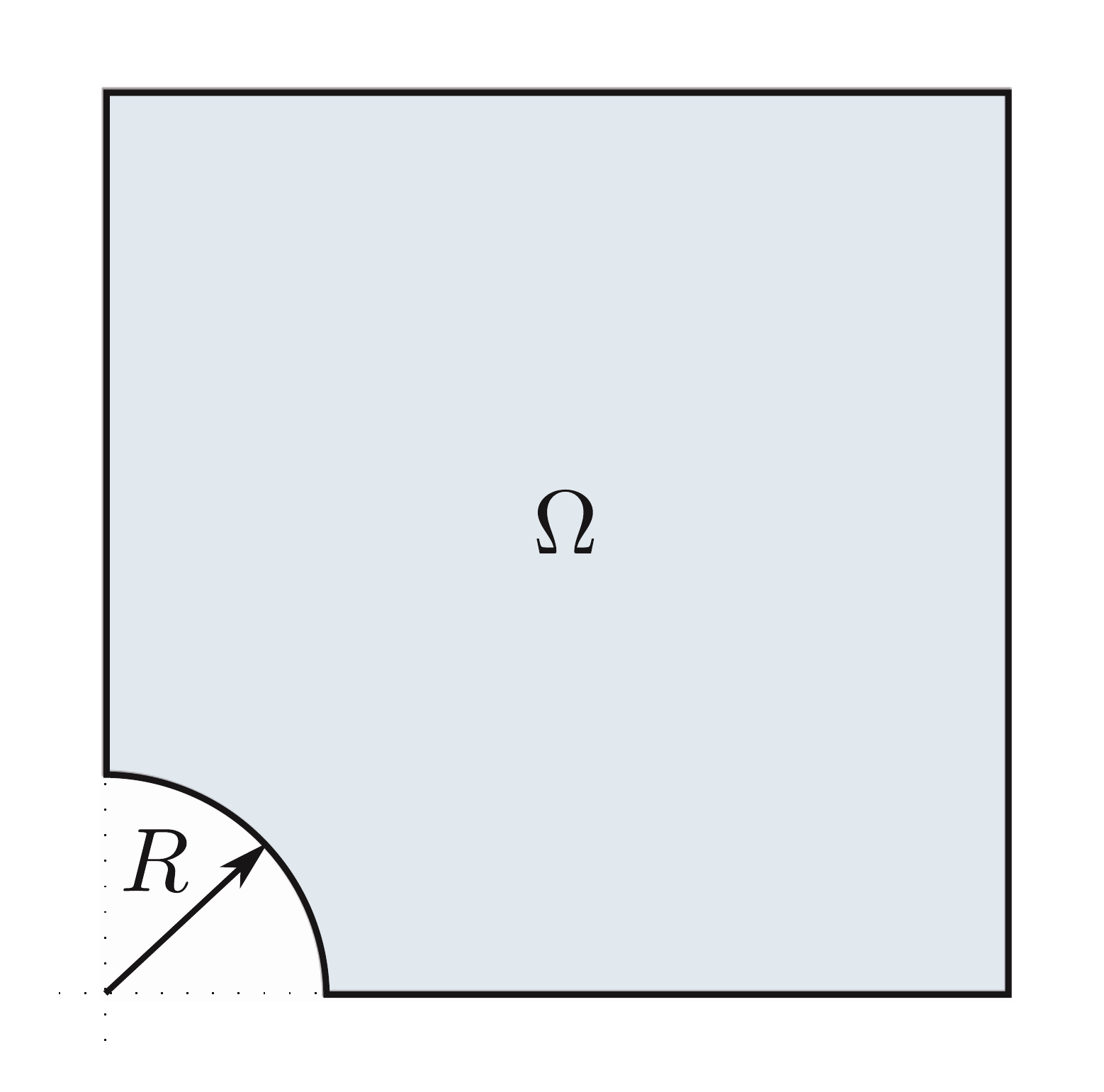}\label{fig:defsPlate}
	\caption{Geometry and boundary conditions of the elastic plate with a hole.  }
	\label{fig:domainDefsPlate}
\end{figure}

%--------------------------------------------------------------------------------
%
As a two-dimensional problem with a non-trivial geometry we compute the infinite elastic plate with a hole subjected to uniaxial tension. The radius of the hole is~$R = 0.25$ and the applied uniaxial traction in the vertical direction is  $ \sigma_\infty = 10^6$.  The Young's modulus of the material is $E = 70 \cdot 10^6$  and its Poisson's ratio is $\nu = 0.3$. This problem has a known closed-form analytic solution~\cite{Timoshenko:1970aa}.  Therefore, we discretise only the plate of size  $L = 1$ shown in Figure~\ref{fig:domainDefsPlate} and apply Dirichlet boundary conditions over its entire boundary. 

\begin{figure} 
	\centering
	\includegraphics[width=0.9\textwidth]{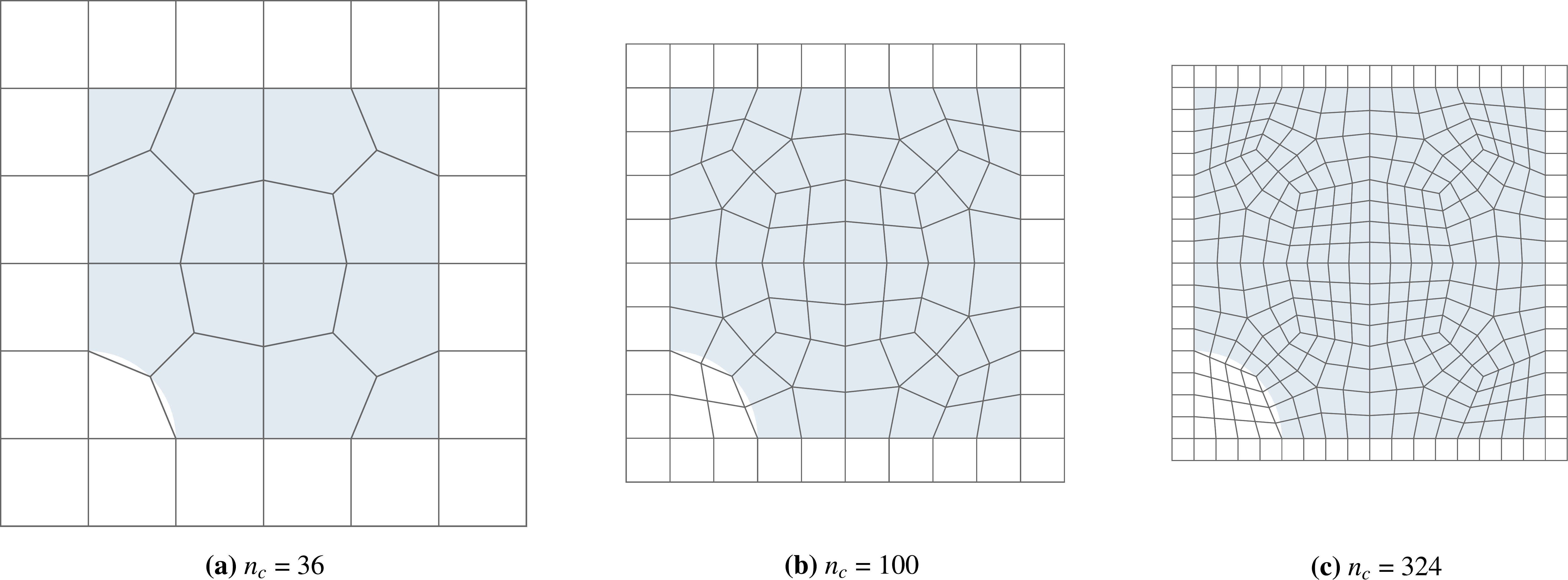}
	\caption{Elastic plate with a hole. Three different partitionings of the plate with a hole (highlighted in blue). The elements outside the domain are the ghost cells. The coarse mesh in (a) is a Voronoi diagram which is refined to obtain (b)  and (c) by introducing new vertices on the cell and edge centres. }
	\label{fig:cellsPlate}
\end{figure}

The initial mesh consists of a Voronoi diagram of 36 non-uniformly distributed points, see Figure~\ref{fig:cellsPlate}. The refined meshes are obtained by subdividing the cells by introducing new vertices on the cell and edge centres. This refinement ensures that the meshes are nested. As can be inferred from Figure~\ref{fig:cellsPlate} along the circular boundary the mesh edges are not aligned with the boundary. In the respective cells cut by the boundary the finite element integrals are evaluated only over the cell areas inside the domain. The cut-cells for integration are obtained with the clipping process introduced in~\ref{app:clipVoro}.  To achieve a higher order approximation the edges of the triangles used for integration are curved by introducing additional nodes on the faces. As in standard finite elements, to achieve an optimally convergent method the boundary geometry has to be approximated with the same polynomial order as the used mollified basis functions.

\begin{figure} 
	\centering
	\includegraphics[width=0.5\textwidth]{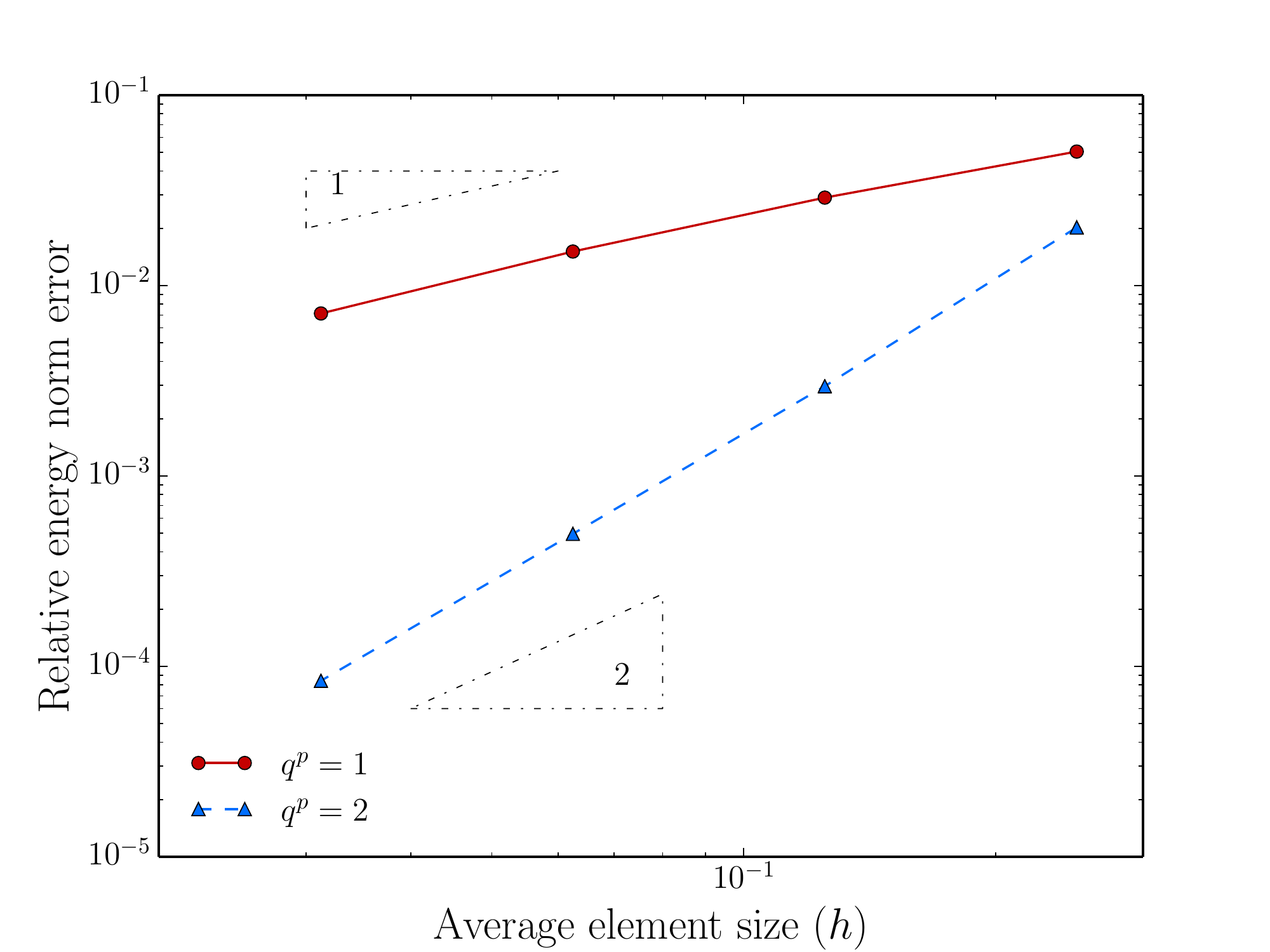}
	\caption{Elastic plate with a hole. Convergence of the relative energy norm error with quartic spline mollifier and local polynomial basis of degree~$q^{p}=1$ and~$q^{p}=2$. }
	\label{fig:convPlate}
\end{figure}

To analyse this problem we again use $C^2$ continuous linear and quadratic basis functions with $q^{p}=1$ and $q^{p}=2$, respectively. Figure~\ref{fig:convPlate} shows the convergence of the errors in the energy norm.  It is apparent that optimal convergence rates are achieved. A final comment concerns the possibly very small contributions of basis functions cut by the boundary to the system matrix. To this end, several approaches have been developed in immersed/embedded finite elements~\cite{Ruberg:2011aa,de2019preconditioning,gurkan2019stabilized}. A particularly simple approach is to scale the relevant basis functions according to their support size within the domain, i.e. 
\begin{equation}
	 \dfrac{ | \supp ( N_i ) \cap \Omega \big ) |}{ |  \supp (N_i) |  } \, .
	\label{eq:fracRatio}
\end{equation}

%

%
%--------------------------------------------------------------------------------          
\subsection{Three-dimensional example}
\begin{figure} 
	%\centering
	\subfloat[Fine surface mesh]{\includegraphics[width=0.4\textwidth]{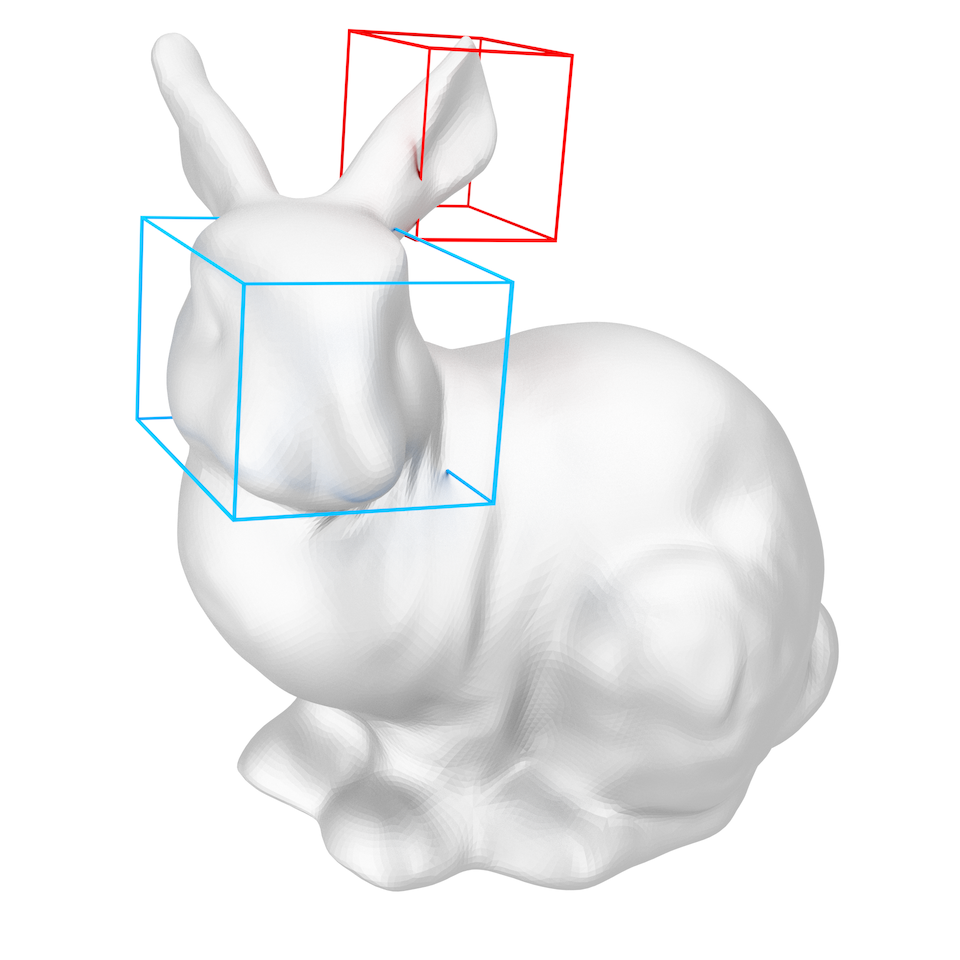} \label{fig:bunnyRefine} }
	\hspace{0.15\textwidth}
	\subfloat[Surface of the clipped Voronoi diagram]{\includegraphics[width=0.3\textwidth]{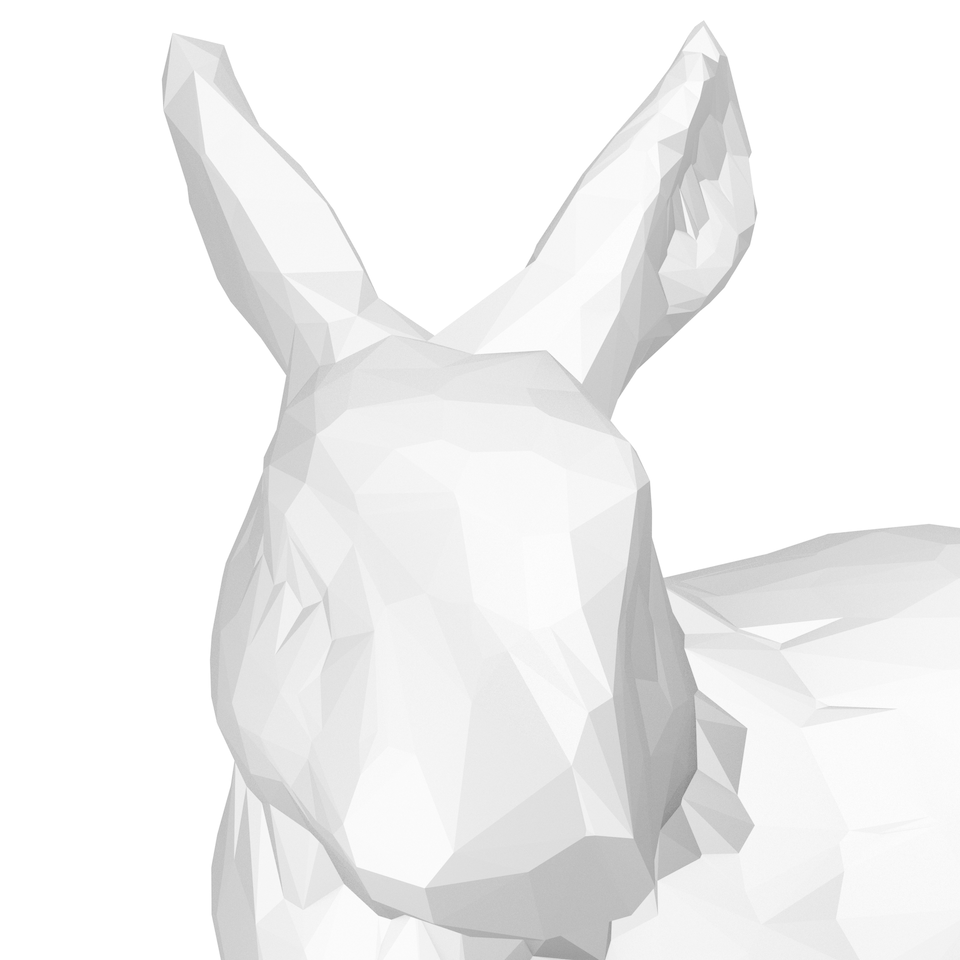} \label{fig:bunnyCloseupRef}}  \\[1cm]
	\centering
	\subfloat[Voronoi diagram]{\includegraphics[width=0.4\textwidth]{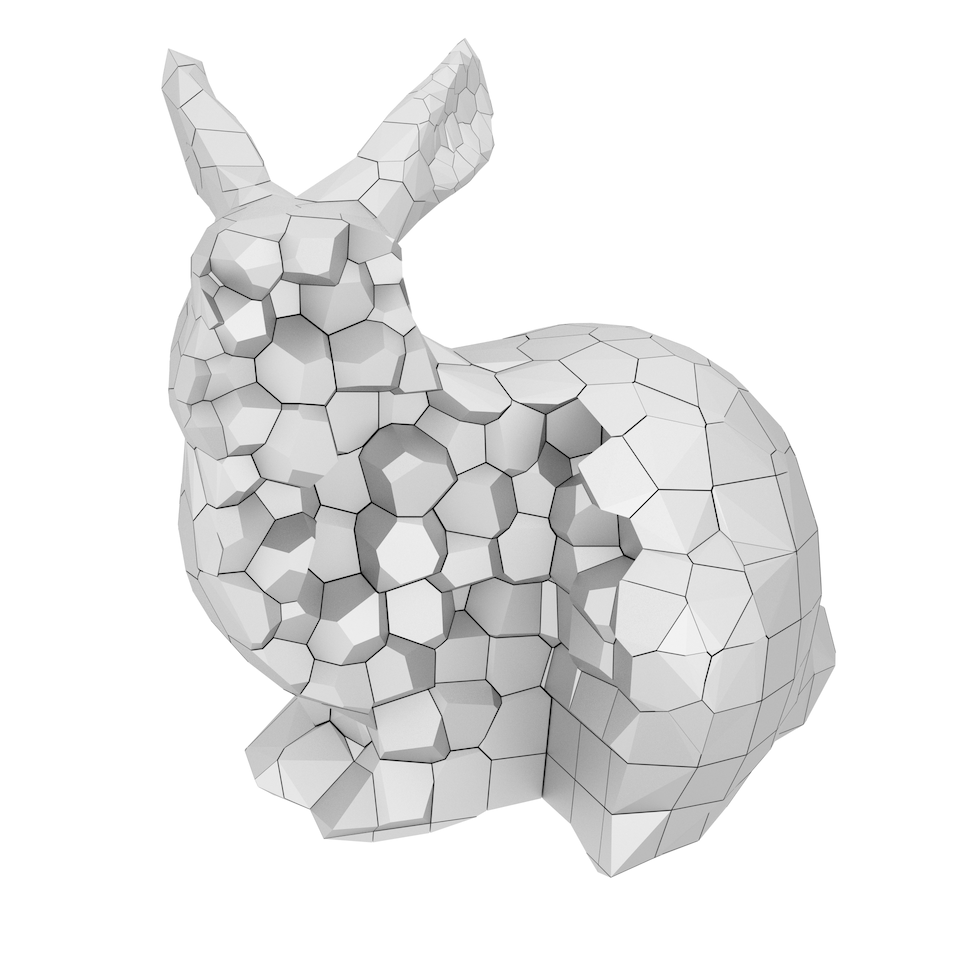} \label{fig:bunnyInternal}}      
	\hspace{0.1\textwidth}
	\subfloat[Finite element solution]{\includegraphics[width=0.4\textwidth]{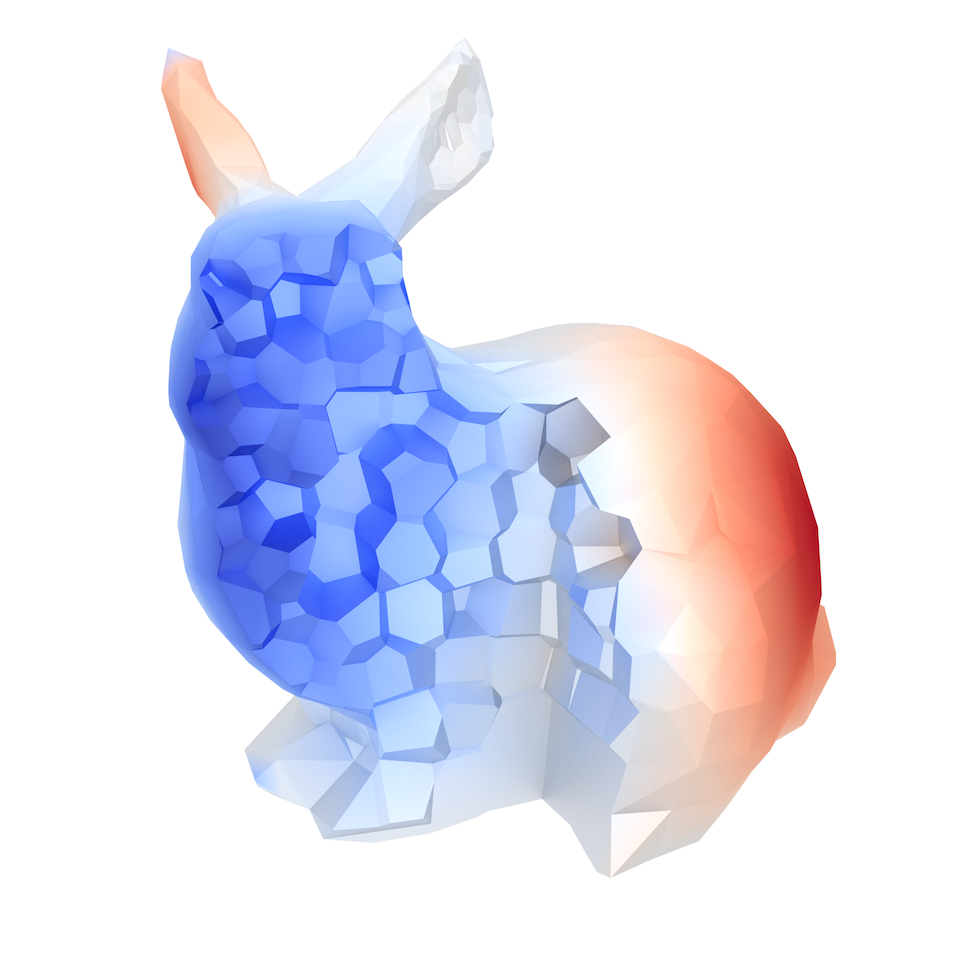} \label{fig:bunnyDisp}}
	\caption{Poisson-Dirichlet problem on the domain contained within the Stanford bunny. The two boxes in (a) indicate the locations where additional nodes are introduced to refine the mesh. The mesh in (c), with some of the cells omitted, shows the unstructured Voronoi diagram used in the mollified finite element computation (d).   \label{fig:bunny}}
\end{figure}

%--------------------------------------------------------------------------------
%
As an illustrative three-dimensional example with a complex boundary we consider the solution of a Poisson-Dirichlet problem on the domain contained within the Stanford bunny, see Figure~\ref{fig:bunny}. The geometry of the bunny is given as a triangle mesh with 66272 facets. The volume mesh shown in Figure~\ref{fig:bunnyInternal} is created in several steps. Firstly we introduce within the bounding box of the bunny a set of uniformly distributed $20 \times 20 \times 20$ points with each $0.08$ apart. The head is then refined by introducing additional $6 \times 7 \times 9$ points with each $0.05$ apart.  Finally, one of the ears is refined by adding $6 \times 7 \times 8$ points with each $0.04$ apart. We then generate the Voronoi diagram of all the points and iteratively relax it to achieve a more even distribution of cell sizes. During the iterative relaxation the new position of each point is recomputed by convolving nodal coordinates with a box function. This relaxation is equivalent to standard Laplace smoothing of meshes. After the relaxation the Voronoi diagram is clipped with the technique described  in \ref{app:clipVoro}. The final mesh consists of 897 cells. 

We solve on the generated polytopic mesh a Poisson-Dirichlet problem with a source term~$s(x)$ such that the solution is equal to 
\begin{equation}
	u = \cos \left( \pi x^{(1)} \right ) \sin \left (\pi x^{(2)} \right ) \cos \left ( \pi x^{(3)} \right )  \, .
\end{equation}
The isocontours of the computed solution are shown in Figure~\ref{fig:bunnyDisp}. For the employed mollified basis functions we use a local polynomial basis with $q^{p}=1$ and a quartic spline mollifier with a support size of~$h_m = 0.16$, which is twice the coarse cell size. The discretisation has in total~$9032$ basis functions. The cut-cells are stabilised as described earlier by scaling the basis functions according to~\eqref{eq:fracRatio}.   

%  
%--------------------------------------------------------------------------------          
\section{Conclusions \label{sec:conclusions}}
%--------------------------------------------------------------------------------
%
We introduced the mollified basis functions of arbitrary order and smoothness and verified their excellent finite element approximation properties with a selected set of examples. In the two- and three-dimensional examples we chose the Voronoi diagram of a given set of points as the partitioning of the domain. The mollified basis functions are obtained by convolving cell-wise defined local polynomial approximants with a compactly supported smooth mollifier with a unit volume. We integrate the convolution integrals exactly (up to round-off errors) by first determining the geometry of the polytopic integration domain and then reducing the dimension of the integrals using the divergence theorem.
%Although we considered only one-, two- and three-dimensional domains, the proposed numerical integration approach to evaluating the convolution integrals applies to any dimension. 
In determining the polytopic integration domain we consider the intersection of a single cell with the mollifier and make use of polytope clipping and convex-hull computation algorithms. Efficient and robust implementations of both algorithms are available in most platforms, including Matlab, Mathematica and Python, and in high-performance geometry libraries~\cite{ray2018meshless}. Because the mollified basis functions are not boundary conforming, we enforce boundary conditions with standard immersed/embedded finite element techniques. The obtained polynomial basis functions may have breakpoints, i.e. points or lines of reduced continuity, within the cells. Therefore, we evaluate the finite element integrals with a variationally consistent approach originally developed for meshless methods. As shown numerically and analytically the mollified basis functions in combination with the proposed finite element implementation can pass the patch test and achieve optimal rates of convergence. Finally, while all the operations in evaluating the mollified basis functions are of geometric nature, in meshless methods, like RKPM, a dense local matrix must be inverted. Evidently, this matrix can become very large especially in 3D when high-order polynomials are used. The algorithmic complexity of the used geometric operations is at most log-linear and the complexity of the matrix inversion is cubic. Therefore, we conjecture that the mollified approximants are more efficient for high-order polynomials and higher dimensions

There are several promising applications and extensions of the proposed mollified approximation scheme worth mentioning. Clearly, it is straightforward to apply h-, p- and hp-refinement. A given Voronoi diagram can be h-refined by incrementally adding new points and updating the Voronoi diagram. For p-refinement it is sufficient to choose in each cell the order of the polynomial approximant differently. A-priori and a-posteriori estimators are crucial for making efficient use of h-, p- and hp-refinement in applications. Furthermore, in our present implementation the mollifier support size is uniform throughout the domain. As our preliminary experiments indicate, it is possible to vary the mollifier support size within a domain. This can be, for instance, used in creating boundary interpolating approximants by continuously shrinking the mollifier support size to zero while approaching the boundary. An alternative approach to easing the enforcement of boundary conditions is to blend mollified basis functions with standard finite elements as in blending techniques developed for meshless methods~\cite{huerta2000enrichment,rosolen2013blending}. Lastly, returning to our original motivation in developing smooth basis functions for isogeometric analysis, it is appealing to develop mollified blending techniques for B-spline patches meeting at an extraordinary vertex. The convolutional definition of the B-splines can be used to derive mollified approximants which reduce to B-splines away from the mesh boundaries.

\appendix
%\newpage

%
%--------------------------------------------------------------------------------          
\section{Convergence estimates \label{app:convergence}}
%--------------------------------------------------------------------------------
%
We make use of convolution  and polynomial approximation estimates to derive convergence estimates for the proposed mollified finite element approximation scheme~\cite{ern2004theory, adams2003sobolev}. In doing so, we adopt a multiindex notation and denote derivatives with~\mbox{$D_j = \partial/\partial x_j$}. If the multiindex \mbox{$\alpha = (\alpha_1,\dotsc, \, \alpha_d)$} is an $d$-tuple of non-negative integers $\alpha_j$, then $D^\alpha = D_1^{\alpha_1} \cdots D_n^{\alpha_d}$ is a differential operator of order $|\alpha| = \alpha_1 + \cdots + \alpha_d$, with the convention that \mbox{$D^{(0,\dots,0)}u = u$}. We denote the norms of the Lebesgue and Sobolev spaces~$L^{p} (\Omega)$ and $W^{p,k}(\Omega)$ with  
\begin{align}
	\begin{split}
		\| f \|_{L^p} & \coloneqq \left (  \int_\Omega | f(\vec x) |^p \D \vec x \right )^{1/p} \\
		\| f \|_{W^{p, k}}  &\coloneqq \left (  \sum_{|\alpha | \le k}  \|  D^{\alpha} f  \|^p_{L^p}   \right )^{1/p}  \qquad 
		| f |_{W^{p, k}}  \coloneqq \left (  \sum_{|\alpha | = k}  \|  D^{\alpha} f  \|^p_{L^p}   \right )^{1/p}  \, .
	\end{split}
\end{align}
In this Appendix we assume that the mollifier is $C^\infty$ continuous. For instance, we may take  
\begin{equation}
    \rho(\vec x)
    =
    \left\{
    \begin{array}{ll}
        \exp[{-1/(1-| \vec x|^2)}]/I_n , & \text{if } |\vec x| < 1 \,   \\
        0, & \text{otherwise} \, ,
    \end{array}
    \right.
\end{equation}
where the scalar $I_n > 0$  is chosen so that~$\rho(\vec x)$ has unit volume. The scaled mollifier with the support size~\mbox{$| \Box_{\vec 0}| = h_m$}  is given by
\begin{equation}\label{JUmus5}
    m (\vec x)  = \rho(2 \vec x/ h_m) 
\end{equation}
so that
\begin{equation}\label{z6Chlw}
    \|D^\alpha m \|_{L^1}
    \leq
    \frac{C}{h_m^{|\alpha|}} 
\end{equation}
for any multiindex $\alpha$.  

As shown, e.g. in~\cite{adams2003sobolev}, the derivatives of the convolution satisfy for any pair of multiindices $\alpha$ and $\beta$ the relation
\begin{equation}\label{br4hEb}
    D^{\alpha+\beta}(m * u) = D^\alpha m * D^\beta u 
\end{equation}
and the Young's inequality for convolutions reads
\begin{equation}
    \|f * u\|_{L^p} \leq  \|f\|_{L^1} \|u\|_{L^p} \, 
\end{equation}
yielding the estimate
\begin{equation}\label{dROd1A}
    \| m * u\|_{L^p} \leq \|u\|_{L^p}  \, .
\end{equation}
Moreover, according to the Bramble-Hilbert lemma, for a  $u\in W^{p,k+1}(\Omega)$, $k\geq 1$, $1\leq p \leq \infty$, there is  a $v_h\in P_k$ (polynomials of degree less than or equal to $k$ in one variable)  such that
\begin{equation}\label{crLB1u}
    \| v_h - u \|_{L^p(\Omega)}
    \leq
    C h^{k+1} |u|_{W^{p,k+1}(\Omega)} \, . 
\end{equation}
Let us now extend $u$ and $v_h$ by zero outside the domain $\Omega$ and denote their mollifications with 
\begin{equation}
   \widehat{ v}_h  = m * v_h  \quad \text{and} \quad \widehat{u} = m  * u  \, .
\end{equation}
Then, for $l < k$, from \eqref{br4hEb}, \eqref{dROd1A}, \eqref{z6Chlw} and \eqref{crLB1u} we have
\begin{equation}
    \| \widehat v_h  - \widehat u  \|_{W^{p,l}(\Omega)}
    \leq
    \frac{C}{h_m^l}
    \| v_h - u \|_{L^p(\Omega)}
    \leq
    \frac{A h^{k+1}}{h_m^l} |u|_{W^{p,k+1}(\Omega)} ,
\end{equation}
where $A>0$ is a constant. The difference between a function and its mollification can be bounded with a standard approximation theorem for convolutions (\cite{adams2003sobolev}, Theorem 5.33) which reads
\begin{equation}
    \| \widehat u  - u \|_{W^{p,l}(\Omega)}
    \leq
    B h_m^{k+1-l}
    |u|_{W^{p,k+1}(\Omega)} ,
\end{equation}
where $B>0$ is a constant. Combining the preceding two estimates we obtain
\begin{equation}\label{1iSplb}
\begin{split}
    \| \widehat v_h  - u \|_{W^{p,l}(\Omega)}
    & \leq
    \| \widehat v_h - \widehat u  \|_{W^{p,l}(\Omega)}
    +
    \| \widehat u - u \|_{W^{p,l}(\Omega)}
    \\ & \leq
    \left( A \frac{h^{k+1}}{h_m^l} + B h_m^{k+1-l}\right)
    |u|_{W^{p,k+1}(\Omega)} .
\end{split}
\end{equation}
This bound is minimised by taking
\begin{equation}\label{1IjUpR}
    h_m
    =
    \left(
        \frac{l A}{(k+1-l) B}
    \right)^{1/(k+1)}
    h
    \equiv h_m^* ,
\end{equation}
which defines the optimal mollifier support size $ h_m^*$ in dependence of the mesh size $h$.  Inserting the optimal~$ h_m^*$ into the estimate~\eqref{1iSplb} gives
\begin{equation}\label{X2bRUd}
    \|  \widehat {v}_h \phantom{}^* - u \|_{W^{p,l}(\Omega)}
    \leq
    C h^{k+1-l} |u|_{W^{p,k+1}(\Omega)} ,
\end{equation}
where~$   \widehat {v}_h \phantom{}^* $ denotes the mollification of the polynomial~$v_h \in P_k$ with a mollifier with support sizel~$h_m^*$.  We note that the estimate \eqref{X2bRUd} provides control over $l$-th order derivatives, whereas the initial estimate \eqref{crLB1u} does not. 

As in standard finite element approximation theory, see e.g.~\cite{ern2004theory},  applying estimate~\eqref{X2bRUd} cell-wise and considering their sum yields global convergence estimates. Subsequently, it is straightforward to confirm the optimal convergence of the proposed mollified finite elements as already suggested by our numerical experiments. 

%
%--------------------------------------------------------------------------------          
\section{Clipped Voronoi diagrams \label{app:clipVoro}}
\begin{figure}[tb] 
	\centering
	\subfloat[Polygonal mesh]{\includegraphics[width=0.3\textwidth]{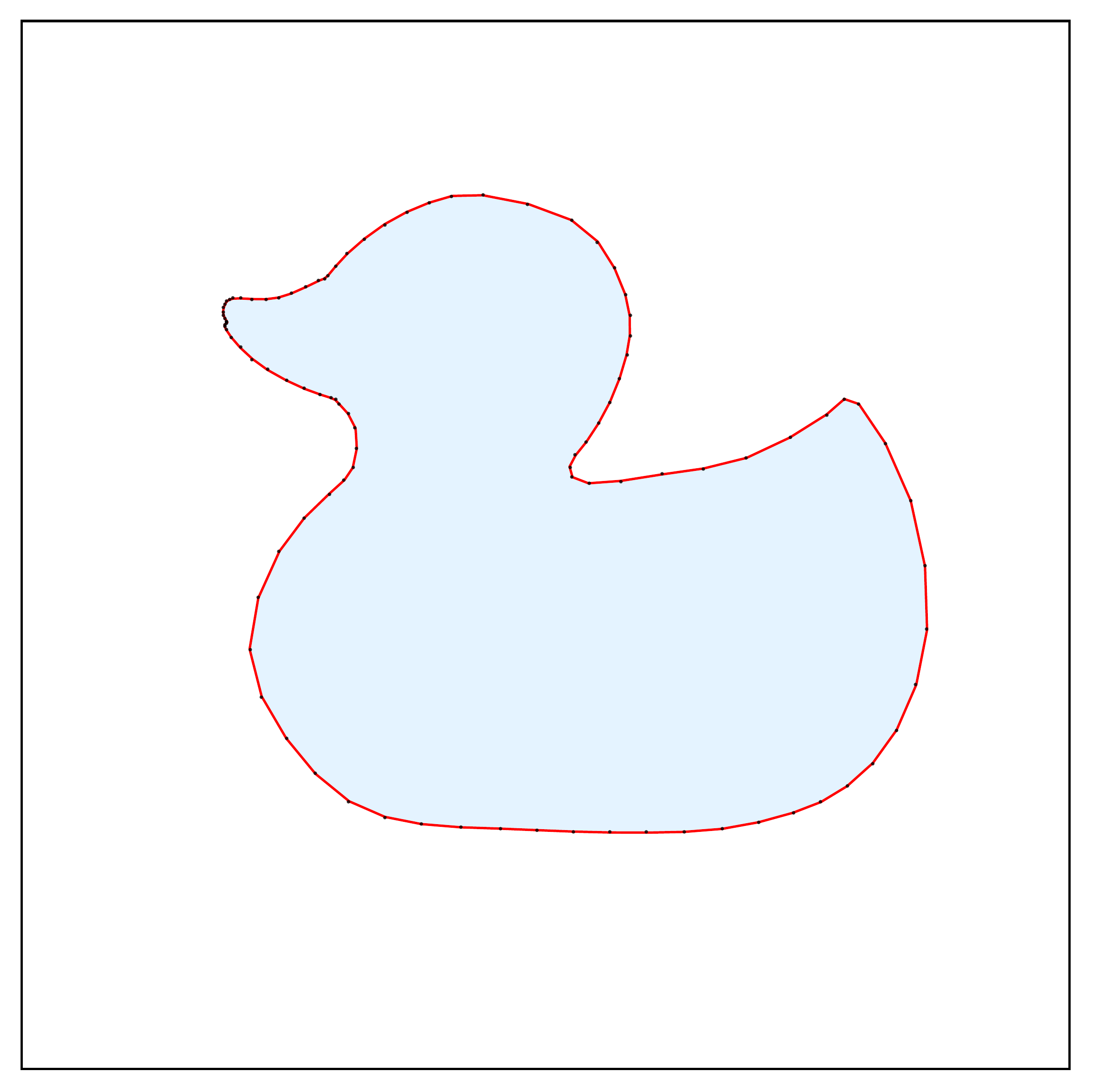}} \hspace{0.1\textwidth}
	\subfloat[Signed distance function]{\includegraphics[width=0.3\textwidth]{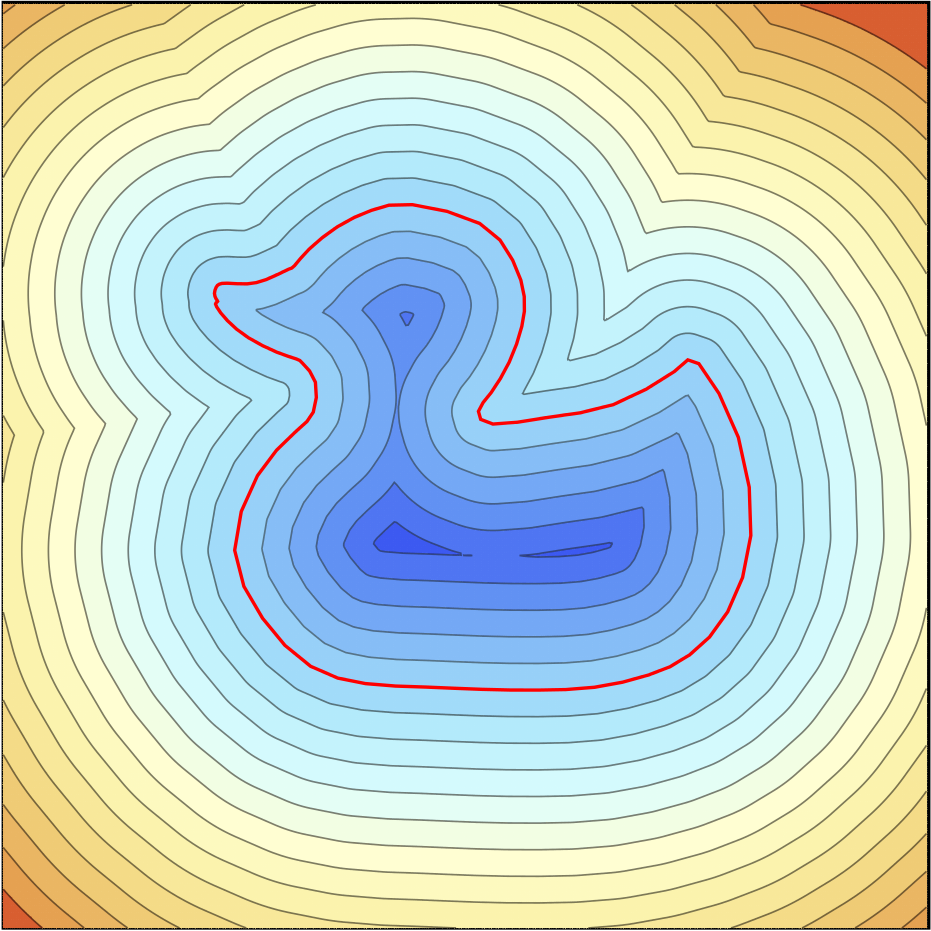}} 
	\\
	\subfloat[Voronoi diagram]{\includegraphics[width=0.3\textwidth]{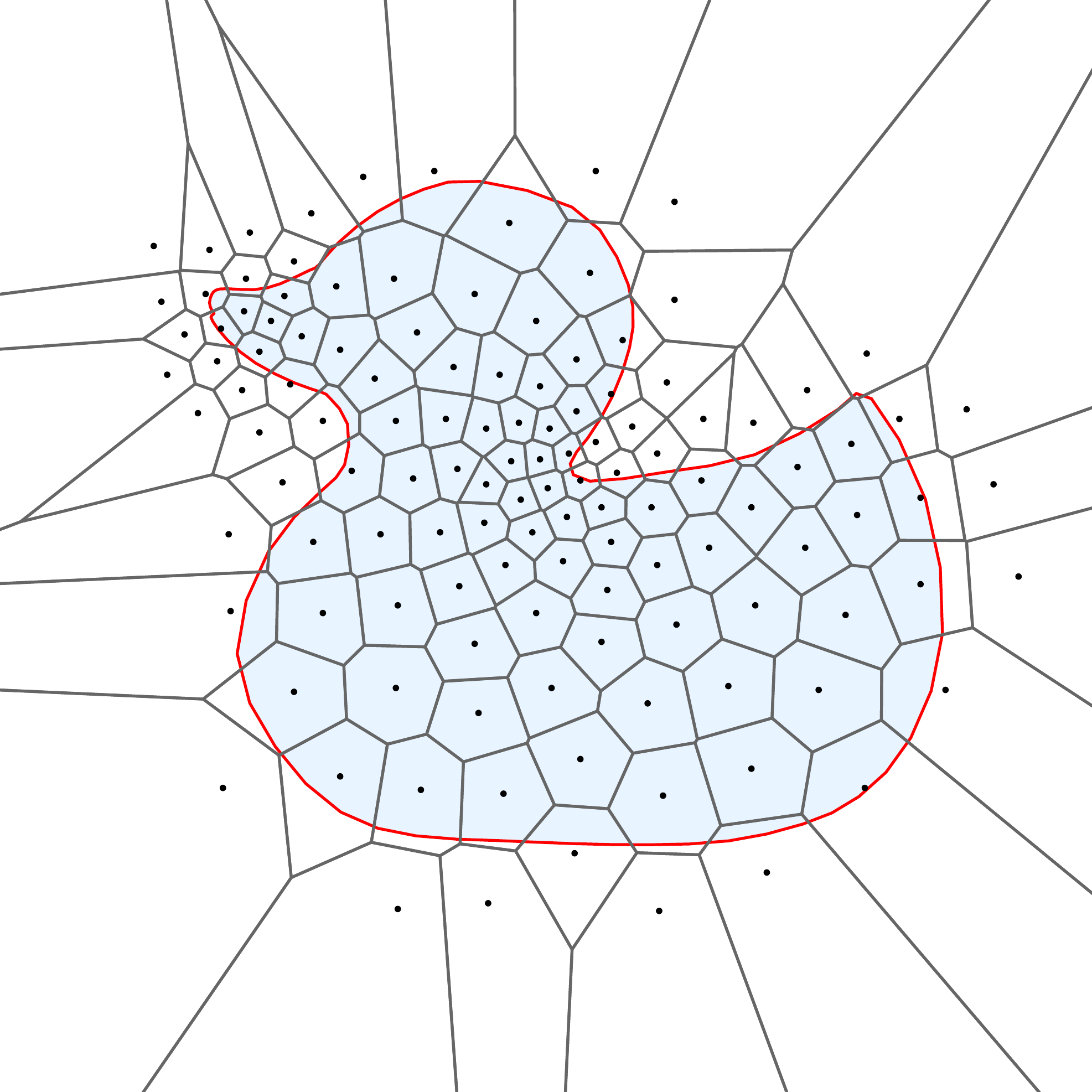}} \hspace{0.1\textwidth}
	\subfloat[Clipped Voronoi diagram]{\includegraphics[width=0.3\textwidth]{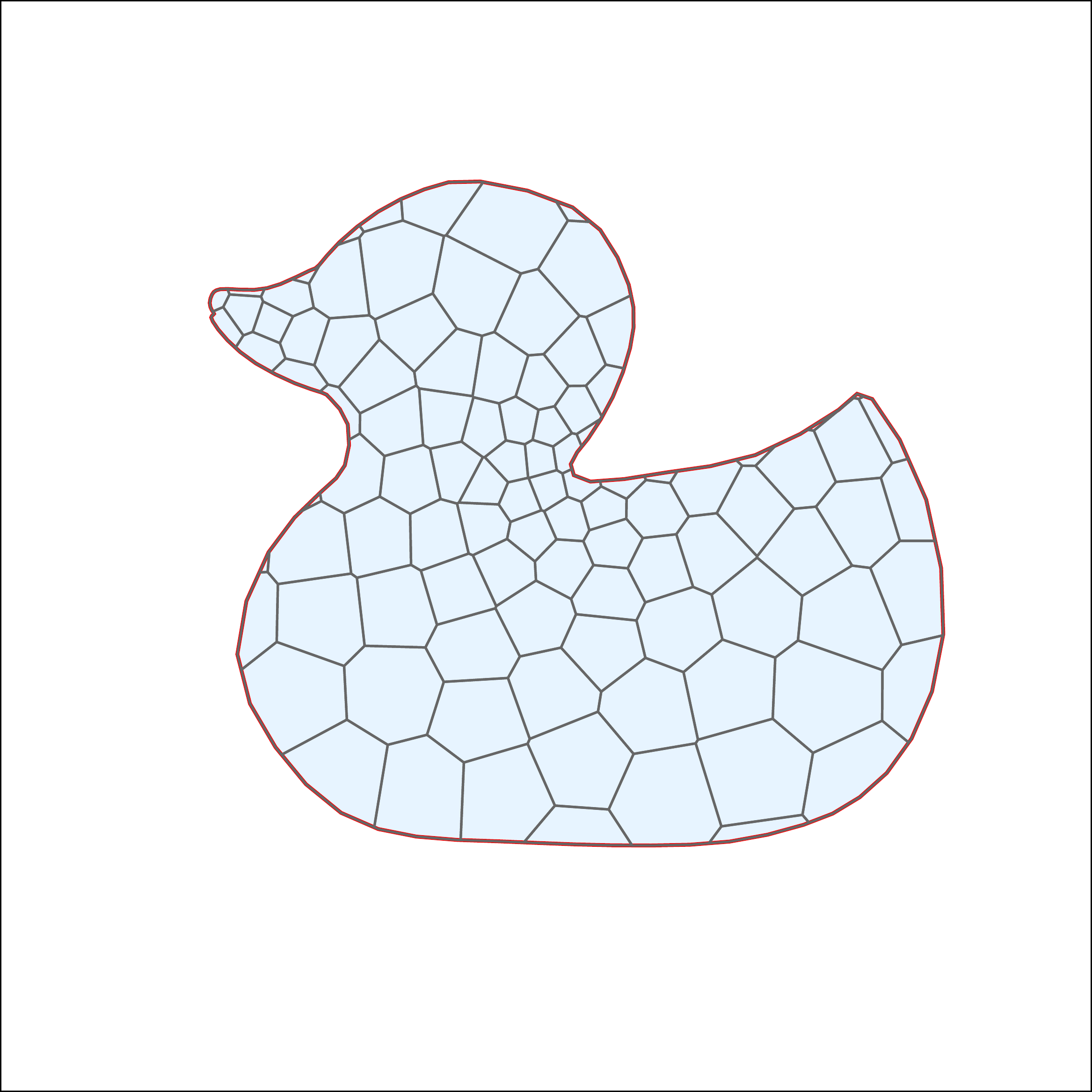}}
	\caption{Clipped Voronoi diagram of a domain described by a polygonal mesh. }
	\label{fig:duckVoronoid2D}
\end{figure}

%--------------------------------------------------------------------------------
%
We briefly review the properties of Voronoi diagrams and sketch the generation of clipped Voronoi diagrams which approximately fill a given domain~$\Omega$. For a more detailed discussion see, e.g.,~\cite{aurenhammer1991, du1999, deBerg2010}. For a set of points   $\{ \vec c_i \}_{i=1}^{n_c}$ in $ \mathbb{R}^d$ the Voronoi diagram is defined by a set of cells $\{ \Omega_i\}_{i=1}^{n_c}$ such that
\begin{equation} \label{eq:vorocrit}
	\Omega_i = \left \{  \vec x \in \mathbb{R}^d | \; |   \vec x - \vec c_i  | \leq | \vec x - \vec c_j |  \quad  \forall i  \neq j  \right \}  \, .
\end{equation}
As indicated in Figure~\ref{fig:duckVoronoid2D}, the cells are convex,  are either bounded or unbounded and have planar faces. There are a number of efficient software libraries available for generating Voronoi diagrams, such as the Voro++~\cite{rycroft2009voro++} library (for 3D) and Mathematica (for 2D) used in this work.  To obtain a Voronoi diagram that approximately fills a given domain~$\Omega$ the cells intersected by the boundary are clipped. To implement the clipping process we assume that the domain~$\Omega$ is described implicitly with a signed distance function
\begin{equation}\label{eq:sgndist}
   \phi (\vec{x}) = 
   \begin{cases} 
     \phantom{-}  	\min_{\vec y \in \Gamma} | \vec x - \vec y | \quad & \text{if } \vec x  \in \Omega \\ 
     \phantom{-}	0 & \text {if } \vec{x} \in \Gamma \\
     				-\min_{\vec y \in \Gamma} | \vec x - \vec y | & \text{otherwise}\, , 
   \end{cases}
\end{equation} 
where~$\Gamma$ is the boundary of the domain~$\Omega$. Domains that are described with a parametric polygonal mesh can first be converted to an implicit signed distance function representation using standard algorithms, see e.g.~\cite{Ruberg:2011aa}.

The minimal data structure for representing a  Voronoi cell $\Omega_i $ consists of its vertices $\{ \vec \eta_{i,j} \}$ and orientated faces $\{ \gamma_{i,j} \}$, with all the face normals~$\vec n_{i,j} $  pointing, e.g., outside the cell. The cells cut by the boundary are determined by evaluating the level set function at the vertices and checking whether 
\begin{equation}
	\min_{j}  \phi (\vec \eta_{i,j}) \cdot \max_{j} \phi( \vec  \eta_{i,j}) < 0 
\end{equation}
is satisfied. Each of the cut cells is clipped by performing the following steps:
\begin{enumerate}
	\item Deduce the set of  cell edges from  $\{ \vec \eta_{i,j} \}$ and $\{ \gamma_{i,j} \}$.
	\item Determine the intersection points between edges and the boundary $\phi(\vec x) =0$ with a bisectioning algorithm. 
	\item Introduce new vertices at the points determined in step 2.
	\item Generate a new clipped cell by determining the convex hull of vertices inside the domain  and the vertices on the domain boundary. 
\end{enumerate}

In finite element computations the clipped cells represent the integration domain for element integrals for the cells crossing the domain boundary. The clipping process introduced yields only convex clipped cells with planar boundaries. For domains with curved boundaries this limits the overall accuracy of finite element method to first  order even when higher order mollified basis functions are used. In this setting, a standard approach to achieving higher order accuracy in immersed finite elements is to curve the planar faces by introducing additional vertices on the faces, which is clearly also applicable to mollified finite elements.   

%
%--------------------------------------------------------------------------------          
\section{Intersection of a convex cell with a box \label{app:intersection}}
%--------------------------------------------------------------------------------
%
The intersection between a convex cell and a box is required to evaluate the mollified basis functions. The box represents the support of the mollifier and is centred at the given evaluation point. The intersection of two convex solids is frequently required in computer graphics and a range of efficient and robust algorithms is available. Our specific approach is motivated by the more general algorithm presented in~\cite{muller1978finding}.  The key idea is that the intersection of a cell and a box can be determined by  clipping the cell in turn by the six half-spaces defining the box. Each half-space is defined by a plane, or its respective normal,  and has an inside and outside. Furthermore, recall that the intersection of two convex solids is always convex. 

As mentioned in~\ref{app:clipVoro}, in our implementation, a cell~$\Omega_i$ is represented by its vertices  $\{ \vec \eta_{i,j} \}$ and orientated faces $\{ \gamma_{i,j} \}$. With this in mind, the sequence of steps in computing the convex polytope representing the intersection domain is as follows:
\begin{enumerate}
	\item Deduce the set of  cell edges from  $\{ \vec \eta_{i,j} \}$ and $\{ \gamma_{i,j} \}$.
	\item Determine the intersection points between the edges and the six half planes in turn while keeping track of the vertices inside the half-spaces.
	\item Generate a new polyhedron by determining the convex hull of vertices inside the domain  and the intersection points on the edges. 
\end{enumerate}

In meshes with a large number of cells it is usually more efficient first to identify the small set of cells which are possibly intersected by a given box. Subsequent intersection computations have to be applied only to the few identified cells. The relevant cells can be efficiently identified with a standard hierarchical bounding volume tree, see e.g.~\cite{klosowski1998efficient}. 

%
%--------------------------------------------------------------------------------          
\section{Minkowski sum of two polytopes \label{app:minkowski}}
\begin{figure}[] 
	\centering
	\subfloat[$\Box_{\vec  0} $ and $\Omega_i$]{\includegraphics[width=0.45\textwidth]{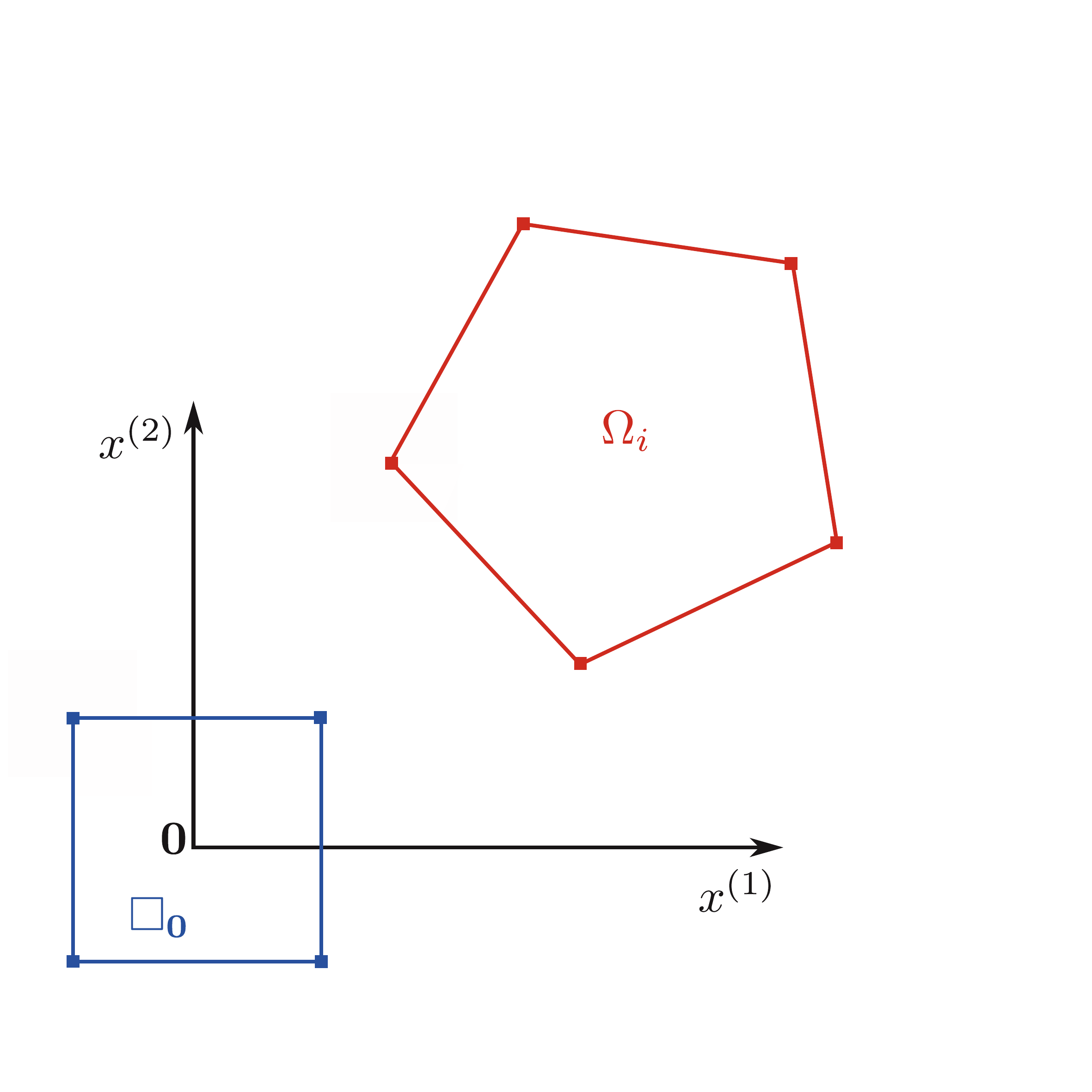}} 
	\vspace{0.05\textwidth}
	\subfloat[$\widehat \Omega = \Box_{\vec  0} \oplus \Omega_i $]{\includegraphics[width=0.45\textwidth]{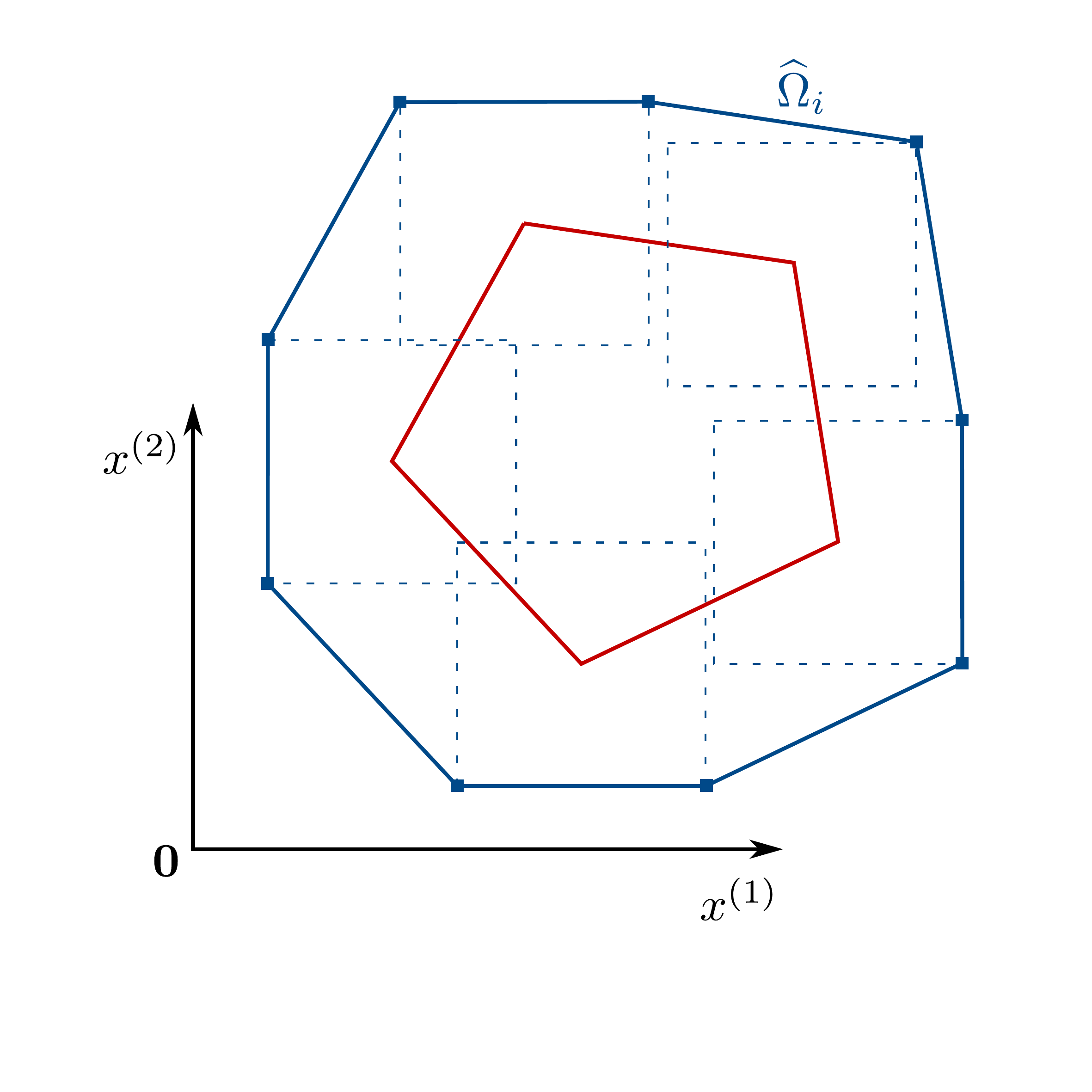} } \\
	\caption{Minkowski sum of  the mollifier support~$\Box_{\vec  0}$ and the cell~$\Omega_i$. The dashed~$\Box_{\vec  \eta_{i, j}}$ in (b) indicate the sliding of the centre of~$\Box_{\vec  0} $ along the boundaries of~$\Omega_i$. }
	\label{fig:minkowski}
\end{figure}

%--------------------------------------------------------------------------------
%
The support of the mollified basis functions~$\vec N_i$ corresponding to the cell~$\Omega_i$ is obtained as the Minkowski sum of the cell with the support of the mollifier.  For an in-depth introduction to Minkowski sums see, e.g.,~\cite{deBerg2010}. The Minkowski sum of two sets~$\Omega_i, \, \Box_{\vec  0} \, \in \mathbb{R}^d$ is defined by
\begin{equation}\label{eq:minkowski}
 \widehat \Omega_i = \Omega_i  \oplus \Box_{\vec  0} = \{ \vec x + \vec y | \,  \vec x \in \Omega_i , \vec y \in \Box_{\vec  0} \} \, .
\end{equation}
The domain~$ \widehat \Omega_i$  resulting from the Minkowski sum may be visualised as that obtained by sliding the centre of~$\Box_{\vec  0} $ along the boundaries of~$\Omega_i$, see Figure~\ref{fig:minkowski}. Recall here that the domain~$\Box_{\vec  0} $ is centred at the origin of the coordinate axis as implied by the subscript~$\vec 0$. It is easy to show that~$ \widehat \Omega_i$ is convex because both ~$\Omega_i$  and~$\Box_{\vec  0}$ are convex. We use this to devise a simple algorithm for computing the Minkowski sum. That is, we first generate a set of points by sliding the domain~$\Box_{\vec  0} $ along the boundaries of~$\Omega_i$, which we subsequently combine  with a convex hull algorithm to obtain~$ \widehat \Omega_i$. In generating the set of points it is sufficient to place~$\Box_{\vec  \eta_{i, j}} $ at the vertices~$\vec \eta_{i,j}$  of the domain~$\Omega_i$ and to take successively the union of the vertices of~$\Box_{\vec  \eta_{i, j}} $, see Figure~\ref{fig:minkowski} . 

\bibliographystyle{elsarticle-num-names}
\bibliography{mollified}

\end{document}